\documentclass[11pt,noamsfonts]{amsart}

\usepackage{amsmath,amsfonts,amssymb,amscd,verbatim,xypic}

\newcommand{\marginlabel}[1]%
  {\mbox{}\marginpar{\raggedleft\hspace{0pt}\bfseries\sf#1}}

\def\ZZ{{\mathbb Z}}

\def\CC{{\mathbb C}}

\def\QQ{{\mathbb Q}}
\def\PP{{\mathbb P}}

\def\cC{\mathcal{C}}
\def\cI{\mathcal{I}}
\def\cJ{\mathcal{J}}

\def\cF{\mathcal{F}}
\def\cG{\mathcal{G}}
\def\cH{\mathcal{H}}

\def\cL{\mathcal{L}}
\def\cO{\mathcal{O}}
\def\cM{\mathcal{M}}
\def\cN{\mathcal{N}}
\def\cP{\mathcal{P}}

\DeclareMathOperator{\Bl}{Bl}
\DeclareMathOperator{\Coker}{coker}
\DeclareMathOperator{\Ker}{ker}

\DeclareMathOperator{\Proj}{Proj}

\DeclareMathOperator{\Spec}{Spec}

\DeclareMathOperator{\codim}{codim}

\newtheorem{lemma}{Lemma}[section]
\newtheorem{theorem}[lemma]{Theorem}
\newtheorem{corollary}[lemma]{Corollary}
\newtheorem{proposition}[lemma]{Proposition}

\theoremstyle{definition}
\newtheorem{definition}[lemma]{Definition}
\newtheorem{example}[lemma]{Example}

\newtheorem{remark}[lemma]{Remark}

\newtheorem{notation}{Notation}
\numberwithin{equation}{section}

\newcommand{\bean}{\begin{eqnarray}}
\newcommand{\eean}{\end{eqnarray}}
\newcommand{\be}{\begin{displaymath}}
\newcommand{\ee}{\end{displaymath}}
\newcommand{\bea}{\begin{eqnarray*}}
\newcommand{\eea}{\end{eqnarray*}}

\newcommand{\ol}{\overline}

\begin{document}

\title{Intermediate Moduli Spaces of Stable Maps}

\author[Andrei Musta\c{t}\v{a}]{Andrei~Musta\c{t}\v{a}}
\author[Magdalena Anca Musta\c{t}\v{a}]{Magdalena~Anca~Musta\c{t}\v{a}}
\address{Department of Mathematics, University of British Columbia,
Vancouver, B.C., Canada}
\email{{\tt amustata@math.ubc.ca, dmustata@math.ubc.ca}}

\date{\today}

\begin{abstract}
 We describe the Chow ring
    with rational coefficients of $\ol{M}_{0,1}(\PP^n,d)$ as
the subring of invariants of a ring
    $B^*(\ol{M}_{0,1}(\PP^n,d);\QQ) $, relative to the action of the
    group of symmetries $S_d$. We compute
   $B^*(\ol{M}_{0,1}(\PP^n,d);\QQ)$ by following
 a sequence of intermediate spaces for $\ol{M}_{0,1}(\PP^n,d)$.
\end{abstract}
\maketitle

\bigskip

\section*{Introduction}

     The moduli spaces of stable maps from curves to smooth
     projective varieties were introduced by
     M. Kontsevich and Y. Manin in \cite{manin}. They provided the
     set-up for an axiomatic algebro-geometric approach to
     Gromov-Witten theory, generating beautiful
     results in enumerative geometry and mirror symmetry.
     Gromov-Witten invariants, defined as intersection numbers
     on the moduli spaces of stable maps, were computed by
     recurrence methods. An important role in these methods was played
     by the ``boundary divisors'' of the moduli space, parametrizing
     maps with reducible domains.

     In the case when the domain curve is rational and the target is $\PP^n $, the functor  $\ol{\cM}_{0,m}(\PP^n , d)$ is represented by a smooth Deligne-Mumford stack,
     $\ol{M}_{0,m} (\PP^n , d).$ Here the generic member is
     a smooth, degree $d$, rational curve in $\PP^n$ with $m$ distinct
     marked points. The boundary is made of degree $d$ morphisms $ \mu : \cC \to
     \PP^n$ from nodal $m$--pointed curves $\cC$ of arithmetic genus
     0, such that every contracted component of $\cC$ has at least 3 special
     points: some of the $m$ marked points or nodes.

     The cohomology ring of $\ol{M}_{0,m} (\PP^n , d)$ is not known
     in general. K.Behrend and A.O'Halloran in \cite{behrend2} have
     outlined an approach for computing the cohomology ring for
     $m=0$, relying on a method of Akildiz and Carell.
     They give a complete set of generators and relations for the case
     $d=2$ and for the ring of  $\ol{M}_{0,0} (\PP^{\infty} , 3)$.

       The main result of this paper is a description of the Chow ring
    with rational coefficients of $\ol{M}_{0,1}(\PP^n,d)$. Our
    method is different from the one employed by K.Behrend and
    A.O'Halloran, relying on a sequence of intermediate moduli spaces. In Theorem 3.23, $A^*(\ol{M}_{0,1}(\PP^n,d);\QQ)$ is expressed as
    the subring of invariants of a
    ring $B^*(\ol{M}_{0,1}(\PP^n,d);\QQ) $, relative to the action of the group of symmetries $S_d$. We
    give a complete set of generators and relations for
    $B^*(\ol{M}_{0,1}(\PP^n,d);\QQ) $, the geometric significance of
    which will be explained  here in more detail.

     Motivated by results in mirror symmetry,
 Givental  in \cite{givental}, and  Lian, Liu and Yau  in \cite{yau} have
 computed Gromov--Witten invariants for hypersurfaces in $\PP^n$,
 via the Bott residue formula for a birational morphism
 $$\Phi: \ol{M}_{0,0}((\PP^n\times \PP^1), (d,1)) \to \PP^n_d.$$
 Here $\PP^n_d:=\PP^{(n+1)(d+1)-1}$ parametrizes $(n+1)$ degree $d$--
 polynomials in one variable, modulo multiplication by constants.
 Following Givental, we will call the domain of $\Phi$ the graph
space. We will use the short notation $G(\PP^n, d)$ for it.

 Of the various boundary divisors of $G(\PP^n, d)$ and their images in $\PP^n_d$, the most notable for
 us is $\ol{M}_{0,1}(\PP^n,d)\times\PP^1$, mapped by $\Phi$ into
 $\PP^n\times\PP^1$. The product $\ol{M}_{0,1}(\PP^n,d)\times\PP^1$
 is embedded in $G(\PP^n, d)$ as the space parametrizing split
   curves $\cC_1\cup\cC_2$, were $\cC_1$ comes with a degree $(d,0)$ morphism
   to $\PP^n\times\PP^1$ and $\cC_2$ comes with a degree $(0,1)$ morphism.
  Our study starts from the diagram
 $$\diagram
 \ol{M}_{0,1}(\PP^n,d)\times\PP^1 \dto \rto & G(\PP^n, d)\dto^{\Phi} \\
 \PP^n\times\PP^1 \rto &
 \PP^n_d.
\enddiagram
$$

A sequence of intermediate moduli spaces $G(\PP^n,d,k)$ is constructed such that the morphism $\Phi:
 G(\PP^n,d)=:G(\PP^n,d,d)\to \PP^n_d=:G(\PP^n,d,0)$ factors through $
 G(\PP^n,d,k+1)\to G(\PP^n,d,k)$. This also induces intermediate spaces 
$\ol{M}_{0,1}(\PP^n,d,k)$. They are described via local $(S_d)^{n+1}$--covers $G(\PP^n,d,k,\bar{t})$
for all homogeneous coordinate systems $\bar{t}$ on $\PP^n$. These are modeled after the construction by
  W.Fulton and R.Pandharipande of the rigidified moduli spaces
  $\left(\ol{M}_{0,0}(\PP^n,d,\bar{t})\right)_{\bar{t}}$, which
  form an \'{e}tale cover  for the stack $ \ol{M}_{0,0} ( \PP^n , d)$.
 
There is a commutative diagram
$$\diagram
 \ol{M}_{0,1}(\PP^n,d,\bar{t})\times\PP^1 \dto \rto & G(\PP^n,
d,\bar{t})  \dto^{\Phi(t)}\rto & \PP^1[(n+1)d] \dto \\
(\PP^1)^n\times\PP^1  \rto &
 \PP^n_d(\bar{t}) \rto & (\PP^1)^{(n+1)d}.
\enddiagram
$$
where $\PP^1[(n+1)d]$ is 
  the Fulton--MacPherson configuration space of $\PP^1$ (see
  \cite{fulton2}). 
 $\PP^n_d(\bar{t})$ is a torus bundle over $(\PP^1)^{(n+1)d}$ and $
 G(\PP^n, d,\bar{t})$ is its pullback to an open subset of the
 configuration space $\PP^1[(n+1)d]$.

  The morphism  $\PP^1[(n+1)d] \to (\PP^1)^{(n+1)d}$ has been
  introduced and described as a sequence of blow--ups by W.Fulton and
  R.MacPherson in \cite{fulton2}. We consider a different sequence of
  blow--ups, more symmetric with respect to the $(n+1)d$
  points. If instead of
  $\PP^1[(n+1)d]$ one considers $\ol{M}_{0,n}$, the resulting sequence is already known (see
  \cite{hassett},  \cite{thaddeus}). In our case we describe
  the intermediate spaces explicitly  as moduli spaces, leading,
  via pull--back of the torus bundles $\PP^n_d(\bar{t})$, to $G(\PP^n,d,k, \bar{t})$. 
Thus the morphism $G(\PP^n,d,k+1, \bar{t})\to G(\PP^n,d,k, \bar{t})$ is simply a composition 
of regular blow--ups, and the blow--up loci are all smooth, equidimensional, 
transverse to each other and are mapped into each other by the natural $(S_d)^{n+1}$--action. 

Following \cite{vistoli}, there is a natural Deligne--Mumford stack associated to the 
coarse moduli scheme $G(\PP^n,d,k)$. An \'etale cover of this stack is constructed 
out of quotients of the $\bar{t}$--covers by small groups. The above properties of 
the blow--up loci are preserved at the level of \'etale cover, but the blow--up is weighted. 
Therefore we will call the morphism 
$G(\PP^n,d,k+1, \bar{t})\to G(\PP^n,d,k, \bar{t})$  a
weighted blow--up along a regular local embedding.

 In this paper we restricted our attention to an intersection--theoretical study of
 $\ol{M}_{0,1}(\PP^n,d)$ and its intermediate spaces
 $\ol{M}_{0,1}(\PP^n,d,k)$. The coarse moduli space and its associated Deligne--Mumford 
stack share the same Chow rings. The stack introduced above finely represents the moduli functor, which is shown in \cite{noi2}. 

The first intermediate space is a
 weighted projective fibration over $\PP^n$, the class of its line bundle
 $\cO(1)$ being the cotangent class $\psi$. The unique polynomial
 $P(t)$  having $\psi$ as root in the cohomology of
 $\ol{M}_{0,1}(\PP^n,d,1)$ can be written in terms of the
 $J$--function of $\PP^n$ as follows: $P(t)=t^{-1}J_d^{-1}$. These calculations are rooted
 in a simple form of Atyiah--Bott localization.
 Then the Chow rings of $\ol{M}_{0,1}(\PP^n,d,k)$ and its substrata are computed
 by induction on $k$. In order to do this, an extended Chow ring is defined for a network generated by regular local embeddings, which morally mimics the Chow ring of the \'etale cover of the stack. There are no complete objects at the level of \'etale covers, but all potential generators of the Chow ring descend and glue to complete
 smooth Deligne--Mumford stacks $\ol{M}^k_I$. The extended Chow ring is generated by the classes of these stacks.

 In our view $\ol{M}_{0,1}(\PP^n,d)$ is the
 natural starting point in the study of the Chow ring of
 $\ol{M}_{0,m}(\PP^n,d)$ for any $m$. On one hand,
 the Chow ring of $\ol{M}_{0,0}(\PP^n,d)$ could be computed from that of
 $\ol{M}_{0,1}(\PP^n,d)$ by analogy with the computation of the Chow
 ring of the
 Grassmannian $\mbox{ Grass }(2,n)$ from that of the flag variety
 $\mbox{ Flag }(1,2,n)$, as we
 hope to show in an upcoming paper. On the other hand, the methods
 of this paper can be applied for the computation of the
 Chow ring of any $\ol{M}_{0,m}(\PP^n,d)$ with $m>1$ (see \cite{noi2}).

  The plan of the paper is as follows: section 1 describes the intermediate spaces $G(\PP^n,d,k)$, $\ol{M}_{0,1}(\PP^n,d,k)$ and the morphisms among them. Section 2 contains a modular presentation of their canonical stratifications.  Section 3 is an extended account of the induction steps involved in the
  computation of the Chow ring of $\ol{M}_{0,1}(\PP^n,d)$, and the Appendix contains some calculations leading to the simplified final formula of the Chow ring.

 This paper is based on an idea of the first author. The detailed construction
 of the intermediate moduli spaces and an outline of the method for
 computing the Chow ring of $\ol{M}_{0,1}(\PP^n,d)$ are contained in Andrei Mustata's Ph.D thesis at the University of
 Utah. Subsequent work by the two authors led to the present
 set--up for the construction of the ring $B^*(\ol{M}_{0,1}(\PP^n,d);
 \QQ)$, and the computations that ensued.

   We warmly thank  Aaron Bertram and Kai Behrend for helpful
 discussions and suggestions.

\bigskip

\section{A description of the graph space}

 Enumerating rational curves in a projective variety $V$
 requires as starting point the existence of a compactification for
 the space of smooth rational curves in $V$. When a curve class
 $\beta\in H_2(V,\ZZ)$ is fixed, Kontsevich--Manin's  moduli space of
 stable maps $\ol{M}_{0,0}(V,\beta )$ provides such a suitable
 compactification. A \emph{stable map} is a tuple $(C, \pi, f)$ where $\pi
 : C\to S $ is a nodal curve over a scheme $S$, $f:C\to V$ is a map
 such that $f_*[C]=\beta$, and every geometric fiber of $\pi $ over a point in $S$ has
 only finitely many automorphisms that preserve the map $f$. 

In particular, the \emph{graph space}
 $G(\PP^n,d):=\ol{M}_{0,0}(\PP^n\times\PP^1, (d,1))$ was the background for the first results in Gromov--Witten theory ( see
 \cite{givental}, \cite{yau}, \cite{bertram}). On the other hand, the space of
 parametrized, degree $d$ smooth rational curves in $\PP^n$ admits
 another, naive compactification given by of
 $\PP^n_d=\PP^{(n+1)(d+1)-1}$, the space of $(n+1)$ degree (at most) $d$
 polynomials in one variable, modulo multiplication by a scalar. The existence of a natural
 birational morphism $$\Phi : G(\PP^n,d) \to \PP^n_d$$
was instrumental in the computations mentioned above.

 Our first
observation is that $\Phi $ may be regarded as a natural
transformation of moduli functors, as we shall discuss in Proposition 1.3. Given a stable map $(C, \pi, f)$ with
 parametrization $\mu:C\to S\times \PP^1$,  the transformation $\Phi (S)$ contracts
 the unparametrized components of $C$. The morphism $f:C\to \PP^n$ induces a rational 
map $f_0:S\times\PP^1\to\PP^n$. The points in the base locus of $f_0$ correspond to contracted components of $C$.

 In the above context, an unparametrized component of a rational curve whose removal does not
 disconnect the curve is called a \emph{tail}.
 One can conceive of various ways to factor the morphism $\Phi$. For
 example, contraction of
 unparametrized components of $C$ could proceed by sequentially
 contracting tails of increasing degree. In the following we define
 moduli problems for the intermediate steps we have just proposed:

\begin{definition} Let $k$ be a natural number, $0\leq k\leq d$. Fix a small rational number $\epsilon$ such that $0<\epsilon <1$.

 A $k$--stable family of degree $d$ maps from parametrized rational curves into $\PP^n$ is a tuple $(C, \mu, \cL, e)$, where $\mu:C\to S\times \PP^1$ is a morphism over $S$ of degree 1 on each fiber $C_s$ over   $s\in S$, $\cL$ is a line bundle on $C$, of degree $d$ on each fiber $C_s$, and $e:\cO^{n+1}\to \cL$ is a morphism of fiber bundles such that:
\begin{enumerate}
\item $( \mu^*\cO_{S\times \PP^1}(1)\otimes\omega_{C|S})^{d-k+\epsilon }\otimes \cL$ is relatively ample over $S$,
\item $\cG :=\Coker e$, restricted to each fiber $C_s$, is a skyscraper sheaf, and \begin{itemize}
\item $\dim \cG_p\leq d-k$ for any $p\in C_s$, where $\cG_p$ is the stalk of $\cG$ at $p$;
\item if $0<  \dim \cG_p$ then $p\in C_s$ is a smooth point of $C_s$.
\end{itemize}
\end{enumerate}

\end{definition}

By condition (1), a $k$--stable map from $C$ to $\PP^n$ will have no
tails of degree less than $d-k+1$, while by condition (2), points
with multiplicity at most $d-k$ may be in the base locus of $e$.
These points replace the smaller degree components of the usual
stable maps. In the case $k=d$, the existence of $\epsilon$ makes formula (1) equivalent to the classical stability condition; in all other cases, $\epsilon$ may be taken to be zero. 

As usual, a morphism between two  families $(C,\mu,\cL, e)$ and
$(C',\mu',\cL',e')$ of $k$--stable, degree $d$ pointed maps consists
of a Cartesian diagram \bea \diagram C \rto^{g} \dto^{\mu} &
C'\dto^{\mu'}\\  S \rto & S'
\enddiagram \eea

such that $g^*\cL'\cong \cL$ and $e$ is induced by $e'$ via the isomorphisms  $\cO_{C}^{n+1}\cong g^*\cO_{C'}^{n+1}$ and $g^*\cL'\cong \cL$.

\begin{notation}
Let $G(\PP^n,d,k)$ denote the functor of isomorphism classes of
$k$--stable, degree $d$ families of maps from parametrized rational
curves to $\PP^n$.
\end{notation}
In the future, the same notation will be somewhat abusively employed
for the corresponding coarse moduli scheme as well as for its moduli
stack. Clearly, $G(\PP^n,d,d)=G(\PP^n,d)$ is the graph space, while
$G(\PP^n,d,0)$ is $\PP^n_d$.

A sequence of intermediate moduli functors is induced for the moduli
space of stable maps with one marked point. The parametrized
component is replaced here by the component containing the marked
point:

\begin{definition} Let $k$ be a natural number, $0< k \leq d$. Similarly
we define a $k$--stable family of degree $d$ pointed maps from parametrized rational curves into $\PP^n$ by a tuple $(C, \pi, p_1, \cL, e)$, where $\pi :C\to S$ is a flat family of genus zero nodal curves over $S$, $p_1: S\to C$ is a section in the smooth locus of $\pi$, $\cL$ is a line bundle of degree $d$ on each fiber $C_s$, $s\in S$, and $e:\cO^{n+1}\to \cL$ is a morphism of fiber bundles satisfying condition (2) from above, such that:
\begin{enumerate}
\item  $\omega_{C|S}(p_1)^{d-k+\epsilon }\otimes \cL$ is relatively ample over $S$,
\end{enumerate}
and such that $\dim \cG_{p_1(s)}=0$ for any $s\in S$.

\end{definition}

\begin{notation}
$\ol{M}_{0,1}(\PP^n,d,k)$ will denote the functor of isomorphism
classes of $k$--stable, degree $d$ families of maps from one--pointed
rational curves to $\PP^n$. Here
$\ol{M}_{0,1}(\PP^n,d,d)=\ol{M}_{0,1}(\PP^n,d)$.
\end{notation}

\begin{proposition}

There are natural transformations of functors $$G(\PP^n, d, k+1) \to
G(\PP^n, d, k),$$ for any natural number $k $ such that $0\leq k<d$,
and
$$\ol{M}_{0,1}(\PP^n, d, k+1) \to \ol{M}_{0,1}(\PP^n, d, k)$$ for $k $ such that
$0<k<d$.

\end{proposition}

\begin{proof}
%de schimbat familia  universala
 Let $(\pi: C_{k+1}\to S, \mu_{k+1}, \cL_{k+1}, e_{k+1})$ be a $(k+1)$--stable family
of degree $d$ maps from parametrized rational curves into $\PP^n$. Set $$\cH:=
( \mu_{k+1}^*\cO_{S\times \PP^1}(1)\otimes\omega_{C_{k+1}|S_{k+1}})^{d-k}\otimes \cL_{k+1}.$$
For sufficiently large $N$, we consider $C_{k}=\Proj (\oplus_{m\geq
0} \pi_* \cH^{mN})$. The sections of $\cH^N$ induce a morphism $f_k:
C_{k+1}\to C_{k} $, which collapses the locus along which $\cH$
fails to be relatively ample, and maps its complement birationally
onto its image. Note that $\cH$ fails to be ample precisely along
tails of degree $d-k$. Let $\cO_{C_{k+1}}(D)$ denote the line bundle
associated to that locus. (This is the pullback of a line bundle on the universal family on the smooth Deligne--Mumford stack $G(\PP^n, d, k+1)$ introduced in Proposition 1.11). Then
$\cL_{k+1}\otimes\cO_{C_{k+1}}((d-k)D)$ descends to a line bundle
$\cL_k$ on $C_k$, with $(n+1)$ global sections determined by the
composition:
$$f_k^*\cO_{C_{k}}^{n+1}\to\cO_{C_{k+1}}^{n+1}\to\cL_{k+1}\to\cL_{k+1}\otimes\cO_{C_{k+1}}((d-k)D).$$
For any geometric point $s\in S$ such that the fiber $C_{k+1 s}$
decomposes as $F_1\cup F_2$ with $F_1\cap F_2=\{ q\}$, and such that
$f_k$ collapses $F_2$ and is an isomorphism on $F_1$, the following
relations over $F_1$:
$$\omega_{C_{k+1 s}| F_1} =\omega_{C_{k s}| F_1}\otimes\cO_{F_1}(q),
$$
$$\cG_{k+1 |F_1} =\cG_{k |F_1}\oplus k_q^{d-k},$$
insure that conditions (1) and (2) in Definition 1.1 are satisfied.

\end{proof}
\bigskip

 Turning our attention to the representability of the functors
 defined above, we recall the notion of rigid stable maps, which
 makes the local structure of the moduli spaces accessible via \'etale atlases.

%surjective birational morphisms
In the case of $G(\PP^n,d)$, the rigid structure accurately
transcribes the construction of  \cite{fulton1} (Proposition 3), with
the only modification that non--parametrized rational curves are
replaced by parametrized ones. Here are the main points:

Let $\bar{t}=(t_0:...:t_n)$ stand for a choice of a homogeneous
coordinate system for $\PP^n$.

\begin{definition}
 A $(\bar{t})$--rigid family of degree $(d,1)$ stable maps from genus 0 curves to $\PP^n\times\PP^1$ consists of data
$ (C, \pi, \mu, \{ q_{i,j} \}_{0 \leq i \leq n, 1 \leq j \leq d} ) $
where
\begin{enumerate}
\item $(\pi \colon \cC \to S , \mu \colon \cC \to \PP^n\times\PP^1)$
is a family of stable, degree (d,1) maps from genus 0 curves to
$\PP^n\times\PP^1$.
\item $( \cC \to S\times\PP^1 , \{ q_{i,j}
\}_{0 \leq i \leq n, 1 \leq j \leq d} )$ is a family of $(n+1)d$
--pointed, genus zero, stable parametrized curves.
\item For any integer $i$ such that $0 \leq
i \leq n$, the following equality of Cartier divisors holds: $$
\mu^*(t_i)=\sum_{j=1}^d q_{i,j}$$
\end{enumerate}
\end{definition}
 \begin{notation}
 The contravariant functor of isomorphism classes of $(\bar{t})$--rigid
families of degree $(d,1)$ stable  maps from rational curves to $ \PP^n \times
\PP^1  $ is denoted by $G(\PP^n,d,\bar{t})$. \end{notation}

\begin{lemma}
$G(\PP^n,d,\bar{t})$ is finely represented by the total space of a
torus bundle over a Zariski open subset of
$\ol{M}_{0,(n+1)d}(\PP^1,1)$.
\end{lemma}

 The fine moduli scheme will also be
denoted by $G(\PP^n,d,\bar{t})$.
 The proof does not differ in any essential way from that of \cite{fulton1},
 Proposition 3. We translate the main steps of the construction to our setup:

  $\ol{M}_{0,(n+1)d}(\PP^1,1) $ is  the Kontsevich--Manin  moduli
space of $(n+1)d$--pointed, degree 1 stable maps and $f:
\ol{M}_{0,(n+1)d+1}(\PP^1,1) \to \ol{M}_{0,(n+1)d}(\PP^1,1) $ is its
universal family, with $(n+1)d$ sections $\{q_{i,j} \}_{0\leq n,
1\leq d}.$

Consider the line bundles $ \mathcal{H}_i =
\mathcal{O}_{\ol{M}_{0,(n+1)d}(\PP^1,1) }(\sum_{j=1}^d q_{i,j} ),$
for $0 \leq i \leq n$. Let $j: B\hookrightarrow \ol{M}_{0,(n+1)d}(\PP^1,1)$
be the open subscheme parametrizing degree 1 pointed maps
$C\to\PP^1$ such that, for every fixed component $C'$ of $C$, the
subsets of the marked point sets $ \{q_{i,j}\}_j $ lying on $C'$
have the same cardinal for all $i\in\{0,...,n\}$. Equivalently,
$ \mathcal{G}_i :=\pi_{B*} \bar{j}^* (\mathcal{H}^{-1}_0 \otimes \mathcal{H}_i )$ is
locally free, and the canonical map $\pi^*_B \pi_{B*} \bar{j}^*
(\mathcal{H}^{-1}_0 \otimes \mathcal{H}_i )\to \bar{j}^*
(\mathcal{H}^{-1}_0 \otimes \mathcal{H}_i )$ is an isomorphism
 for every $1 \leq i \leq n$, where
 $\bar{j}, \pi_B$ are the canonical projections from $ B\times_{\ol{M}_{0,(n+1)d}(\PP^1,1)}\ol{M}_{0,(n+1)d+1}(\PP^1,1)$ to $\ol{M}_{0,(n+1)d+1}(\PP^1,1)$
 and $ B$, respectively.
We say that the map $j$ is balanced. Moreover, every balanced morphism $X \to\ol{M}_{0,(n+1)d}(\PP^1,1)$ factors through $j$. Let then
  $ \tau_i :  Y_i \to B $ denote the total space of the canonical $ \CC^* $ --bundle associated to $ \mathcal{G}_i $.
The torus bundle $ Y := Y_1 \times_B Y_2 \times_B \ldots \times_B Y_n $  represents the functor $G(\PP^n,d,\bar{t})$.

 The $(n+1)$--th power $G:=(S_d)^{n+1}$ of the group of permutations $S_d$ acts on $G(\PP^n,d,\bar{t})$ as permutations of the marked points $\{q_{i,j}\}$. The coarse moduli scheme $G(\PP^n,d)$ is obtained by gluing
quotients of $G(\PP^n,d,\bar{t})$ by $G$, for various homogeneous coordinate systems $\bar{t}$. Concomitantly, $\bigsqcup_{\bar{t}}G(\PP^n,d,\bar{t})$ is an \'etale atlas of the stack $G(\PP^n,d)$.
%explicatii aici?

\bigskip

 There is also an appropriate notion of $\bar{t}$--rigidity for $k$--stable maps, leading to the representability of $G(\PP^n,d,k)$.
The case of $\PP^n_d(\bar{t})=G(\PP^n,d,0,\bar{t})$ is most at hand. The data $(S\times\PP^1, \mu, \{ q_{i,j} \}_{0 \leq i \leq n, 1 \leq j \leq d} ) $
make a $\bar{t}$--rigid, $0$--stable family of maps if conditions (1) and  (3) in Definition 1.4 are satisfied, with 0--stability replacing  the usual stability in (1).
No incidence conditions need be placed on $\{ q_{i,j} \}_{i,j}$. There is a natural, rank $(n+1)$ torus bundle over $(\PP^1)^{(n+1)d}$ representing $\PP^n_d(\bar{t})$, the functor of isomorphism families of $\bar{t}$--rigid, 0--stable maps.
Indeed, an algebraic map $\PP^1\to\PP^n$ is specified by $(n+1)$ independent hyperplane sections and $(n+1)$ constants, the hyperplane sections providing the roots  of $(n+1)$ polynomials, and the $(n+1)$ constants acting as coefficients.
For a more rigorous construction, set  $ P:= ( \PP^1 )^{(n+1)d} $ and for all pairs $(i,j)$ such that $0\leq i\leq n$, $1\leq j\leq d$, let $ \pi_{i,j} : P \to \PP^1 $ be the  $(i,j)$--th projection.
 Define
\[ \cF_i := \pi^*_{i,1} ( \cO_{\PP^1}(1) ) \otimes \pi^*_{i,2} ( \cO_{\PP^1}(1) ) \otimes \cdot \cdot \cdot \otimes \pi^*_{i,d}( \cO_{\PP^1}(1)) \]
and let $ F_i $ be the total space of the canonical $ \CC^* $--bundle
associated to $ \cF_i \otimes \cF_0^{-1} $. Then a standard translation of the arguments in Proposition 3 of \cite{fulton1}, $\PP^n_d(\bar{t})$ is finely represented by $F:=F_1 \times_P F_2
\times_P ...\times_P F_n. $ Notice that if $\pi_P:P\times\PP^1\to P $ is the "universal family" over $P$, and $\{s_{i,j}\}_{i,j}$ are its canonical sections, then $\cO_{P\times\PP^1}(s_{i,j})=\pi_P^*\cF_i\otimes\pi_1^*\cO_{\PP^1}(d)$, hence the analogy with \cite{fulton1}.

\begin{lemma}
The group $G=(S_d)^{(n+1)}$ acts on the universal family of $ \PP_d^n(\bar{t})$ by permutations
  of the sections $\{ q_{i,1}, q_{i,2},..., q_{i,d} \}$,
  for each $0 \leq i \leq n $. Consider the induced action on $ \PP_d^n(\bar{t})$.  The GIT quotients
  $ \PP_d^n (\bar{t}) / G $ for various $ \bar{t} $--s glue together to $ \PP^n_d$, which finely represents the functor
  $ G(\PP^n,d,0) $. \end{lemma}
\begin{proof}
By the reasons spelled out in Proposition 4 of \cite{fulton1}, the quotients $ \PP_d^n (\bar{t}) / G $ for
 various $\bar{t}$ do patch together to a coarse moduli scheme. Moreover,
this is a smooth scheme, $ \PP_d^n (\bar{t}) / G $ being a torus bundle over
$(\PP^d)^{n+1}$. The essential difference here is the existence of a universal
family over the scheme $\PP^n_d$:  a canonical rational  map $\mu:\PP^n_d\times\PP^1\to
\PP^n$. Given any point $(p,x)\in\PP^n_d\times \PP^1$, this map evaluates 
each of the $(n+1)$ degree $d$ polynomials corresponding to $p$ at  the point
$x\in\PP^1$. For a fixed coordinate system  $\bar{t}=(t_0:...:t_{n+1})$, let
$U_{\bar{t}}\subset\PP^n_d$ be the open subset made of points $p$ such that
$\mu(\{p\}\times\PP^1)\not\subset\bigcup_i(t_i=0)$. There is an obvious
bijection $U_{\bar{t}}\to\PP_d^n (\bar{t}) / G $ and, since both spaces are
smooth, an isomorphism. 
\end{proof}
\bigskip

Thus both the graph space $G(\PP^n,d)$ and $\PP^n_d$ admit local finite
covers $G(\PP^n,d,\bar{t})$ and $\PP^n_d(\bar{t})$, respectively. There is a
birational morphism $\Psi:\ol{M}_{0,(n+1)d}(\PP^1,1) \to (\PP^1)^{(n+1)d}$, the
product of all the evaluation maps. Due to the canonical nature of the previous
constructions, the torus bundle $G(\PP^n,d,\bar{t})$ over  $B\subset 
\ol{M}_{0,(n+1)d}(\PP^1,1)$, is the pullback of $\PP^n_d(\bar{t})$ via $\Psi_{|
B}$. 

 We will now see how  the schemes $G(\PP^n,d,k)$ may be obtained by
gluing quotients by $(S_d)^{n+1}$ of torus bundles
$G(\PP^n,d,k,\bar{t})$ over some bases. To find appropriate bases for these
bundles, one needs to factor the morphism $\Psi$ into intermediate steps:

Set $N:=(n+1)d$. The space $\ol{M}_{0,N}(\PP^1,1)$ is also known as the
Fulton--MacPherson compactified   configuration space of $N$ distinct
points in $ \PP^1$, namely the closure (through the diagonal embedding) of  
$(\PP^1)^N\backslash\bigcup_{i,j}\Delta_{i,j}$   in
$(\PP^1)^N\times\prod_{|S|\geq 2}Bl_{\Delta^S}(\PP^1)^S,$ where  
$\Delta_{i,j}$ are the large diagonals in $(\PP^1)^N$, and   for all subsets
$S \subset \{1,...,N\}$, $Bl_{\Delta^S}(\PP^1)^S$ denotes the blow--up of the
corresponding Cartesian product along its small diagonal. From now on this
compactification will be denoted by $ \PP^1[N] $.  Thus $ \PP^1[N] $ may be
constructed from $(\PP^1)^N$ by a sequence of blow--ups of (strict
transforms of) all diagonals in $(\PP^1)^N$. A different order of
blow--ups than the one considered by W.Fulton and R.MacPherson in
\cite{fulton2} also leads to the same result, while the intermediate steps
 are exactly the bases of $\{G(\PP^n,d,k,\bar{t})\}_k$ previously sought:

\begin{notation}
 Let $K$ be an integer such that $0\leq K\leq N$.  $\PP^1[N,K]$ will denote
the  closure through the diagonal embedding of  
$(\PP^1)^N\backslash\bigcup_{i,j}\Delta_{i,j}$   in
$(\PP^1)^N\times\prod_{|S|>N-K}Bl_{\Delta^S}(\PP^1)^S,$ where  
$S$ stands for subsets $S \subset \{1,...,N\}$.  \end{notation}   

The point here is that for all subsets $S$ such that $|S|=N-K$, all the strict
transforms $\tilde{\Delta}^S$ in $\PP^1[N,K]$ intersect transversely, and
$\PP^1[N, K+1]$ is obtained by blowing up $\PP^1[N,K]$  along these
transforms, in any order. Moreover, at each step, the scheme $\PP^1[N,K]$
finely represents a moduli problem of stable parametrized pointed rational
curves:

\begin{definition}
 Consider a morphism $\phi:C\to S\times\PP^1$ of degree 1 over each geometric
fiber $C_s$ with $s\in S$, and $N$
marked sections of $C\to S$. The morphism will be called $K$--stable if
the following conditions are satisfied for all $s\in S$:

 i) not more than $N-K$ of the marked points in $C_s$ coincide;

 ii) any ending irreducible curve in $C_s$, except the parametrized one,
contains more than $N-K$ marked points (here an ending curve is a curve such
that, if removed from $C$, the remaining curve is connected);

 iii) all the marked points are smooth points of the curve $C_s$ and $C_s$ has
finitely many automorphisms preserving the marked points and the map to
$\PP^1$.

\end{definition}

  A similar definition of $K$--stability can be given for unparametrized rational  pointed curves, but then one needs to choose a special point among the marked points; let it be denoted by $P_1$. This marked point is not allowed to coincide with any other marked point. In condition ii), the parametrized component of $C$ is replaced by the component containing $ P_1$.

\begin{proposition}
The smooth scheme $\PP^1[N,K]$ finely represents the functor of isomorphism
families of $K$--stable parametrized rational curves.
\end{proposition}

\begin{proof}
The proposition may be checked by increasing induction on $K$. Given the
universal family $U[N,K]\to\PP^1[N,K]$ with sections $\{s_i\}_{1\leq i\leq
N}, $ then the universal family $U[N,K+1]$ is constructed by blowing up
$\PP^1[N,K+1]\times_{\PP^1[N,K]}U[N,K]$ along the loci where $N-K$ of the
strict transforms $\tilde{s}_i: \PP^1[N,K+1]\to
\PP^1[N,K+1]\times_{\PP^1[N,K]}U[N,K]$  intersect. The transversality
properties of these loci lead to the stability conditions in Definition 1.7.
  
 Given any $(K+1)$--stable family of parametrized rational curves $C\to
S\times\PP^1$, one may contract the tails with only $N-K$ marked points and
obtain a $K$--stable parametrized curve $C'$ over $S$. The fact that the
morphisms $S\to \PP^1[N,k]$ and $C'\to U[N,K]$ canonically lift to $S\to
\PP^1[N, K+1]$ and $C\to U[N,K+1]$ results from the
universality property of the blow--up. 
\end{proof}

For example, $\PP^1[N,0]=(\PP^1)^N$ and $\PP^1[N,N]=\PP^1[N]$.

Let $\bar{t}=(t_0:...:t_n)$ stand for a choice of a homogeneous
coordinate system for $\PP^n$.
We are now ready to define the notion of $\bar{t}$--rigid, $k$--stable family of
maps by relaxing conditions (1) and (2) in Definition 1.4: 

\begin{definition}
A $k$--stable family of degree $d$, $\bar{t}$--rigid maps from
 rational parametrized curves into $\PP^n$ is a tuple $(C, \mu, \cL,
e, \{ q_{i,j}\}_{0 \leq i \leq n, 1 \leq j \leq d})$, such that:
\begin{enumerate}
\item $(C, \mu, \cL, e)$ is a $k$--stable family of degree $d$ maps from
parametrized rational curves into $\PP^n$;

\item  $(C, \mu, \{ q_{i,j} \}_{0 \leq i \leq n, 1 \leq j \leq d} )$ is
an $(n+1)k$ -- stable family of parametrized curves;

\item via the natural isomorphism $H^0(\PP^n, \cO_{\PP^n}(1)) \cong
H^0(\cC,\cO_{\cC}^{n+1})$,
$$ ( e( \bar{t}_i)) = \sum_{j=1}^d q_{i,j}. $$

\end{enumerate}

\end{definition}

 \begin{notation}
 The contravariant functor of isomorphism classes of $(\bar{t})$--rigid
families of degree $(d,1)$, $k$--stable  maps from rational curves to $ \PP^n
\times \PP^1  $ is denoted by $G(\PP^n,d,k,\bar{t})$. \end{notation}

\begin{lemma}
$G(\PP^n,d,k,\bar{t})$ is finely represented by the total space of a
torus bundle over a Zariski open subset of
$\PP^1[(n+1)d, (n+1)k]$.
\end{lemma}

 The fine moduli scheme will also be
denoted by $G(\PP^n,d,k,\bar{t})$.
 The basis of the torus bundle is the open subset $U_k\subset \PP^1[(n+1)d, (n+1)k]$ parametrizing tuples $(C, \{ q_{i,j}\}_{0 \leq i \leq n, 1 \leq j \leq d})$, such that the line bundle $\cO_C(\sum_{j=1}^d q_{i,j})$ does not depend on $i$. $G(\PP^n,d,k,\bar{t})$ is pullback of $\PP^n_d(\bar{t})$ to $U_k$.

\bigskip

 We have discussed how $\ol{M}_{0,1}(\PP^n,d)$ is embedded into $G(\PP^n,d)$.
A rational pointed unparametrized curve, attached to a parametrized $\PP^1$
by identifying the marked point with $0\in\PP^1$ becomes parametrized. 
Accordingly, the spaces $\ol{M}_{0,1}(\PP^n,d,k)$ introduced
in Definition 1.2 also admit local finite $(S_d)^{n+1}$--covers 
$\ol{M}_{0,1}(\PP^n,d,k,\bar{t})$, which are torus bundles over some bases.
$\ol{M}_{0,1}(\PP^n,d, \bar{t})$   for instance, has an open subset in $\ol{M}_{0,(n+1)d+1}\hookrightarrow
\ol{M}_{0,(n+1)d}(\PP^1,1)$ as basis, and $\ol{M}_{0,1}(\PP^n,d,1,\bar{t})$
stands over $\PP^{(n+1)d-2}$, a fiber in the exceptional divisor of
$\PP^1[(n+1)d,1]=\Bl_{\Delta }(\PP^1)^{(n+1)d}$. Definition  1.9 
may be adapted to  $\ol{M}_{0,(n+1)d+1}(k)$, by replacing the parametrized
component with the one containing the first marked point, and asking that this
point is always distinct from the others. In \cite{kapranov} Kapranov has
already described the construction of  $\ol{M}_{0,N+1}$ by successive blow--ups
of $\PP^{N-2}$, while Thaddeus remarked in \cite{thaddeus} that the order of
blow--ups is not important. The intermediate steps employed by us are a
particular case of Hassett's weighted pointed stable curves. (\cite{hassett}).

\begin{proposition}
The coarse moduli schemes $G(\PP^n,d,k)$  are obtained by gluing $(S_d)^{n+1}$--quotients of the rigidified moduli spaces $G(\PP^n,d,k,\bar{t})$, for all homogeneous coordinate systems $\bar{t}$ on $\PP^n$. The same is true for $\ol{M}_{0,1}(\PP^n,d,k)$. There are natural structures of smooth Deligne--Mumford stacks associated to both $G(\PP^n,d,k)$ and $\ol{M}_{0,1}(\PP^n,d,k)$.
\end{proposition}

\begin{proof}
 We have discussed how the schemes $G(\PP^n,d,k)$ and $\ol{M}_{0,1}(\PP^n,d,k)$ are locally quotients of smooth  varieties by the finite group $(S_d)^{n+1}$. Consequently, they are the moduli spaces
of some smooth Deligne--Mumford  stacks (Proposition 2.8 from \cite{vistoli}).
The same notations will be used for these stacks as for their schemes.  
Chevalley--Shepard--Todd  in \cite{todd} show that the quotient of  a smooth
variety by a small group (a group generated by elements that  are invariant
along some divisors) is a smooth variety. Thus for general $k$, an \'etale
cover of the stack $G(\PP^n,d,k)$ will be obtained by factoring an appropriate
neighborhood of any point $x\in G(\PP^n,d,k,\bar{t})$ by the 
largest small normal subgroup $H_x$ of the stabilizer $Stab_x \subset
(S_d)^{n+1}$. This is exactly the subgroup whose elements fix the points of the universal curve over $x$, and thus fail to contribute to the stacky structure of the space of non--rigid $k$--stable maps. 
A detailed proof that $G(\PP^n,d,k)$ and $\ol{M}_{0,1}(\PP^n,d,k)$ are the stacks of $k$-- stable maps is written in \cite{noi2}.

 When $k=d$, there are no small subgroups in $Stab_x$, and the
various $G(\PP^n,d,\bar{t})$ form an \'etale cover of the stack $G(\PP^n,d)$.
For $\PP^n_d(\bar{t})$, the entire group $(S_d)^{n+1}$ is small and $\PP^n_d$
is obviously smooth. For $k$ in the middle, the map from
$G(\PP^n,d,k,\bar{t})$ to $G(\PP^n,d,k)$  factors as a composition of a $GIT$
quotient with an \'etale map of stacks.
\end{proof}

\bigskip

\section{Substrata of intermediate spaces}

\begin{notation}
From now on, the moduli space $\ol{M}_{0,1}(\PP^n,d,k)$ will be denoted simply by $\ol{M}^k$.
\end{notation}

There is a canonical stratification of $\ol{M}_{0,1}(\PP^n,d)$ (and similarly of the graph space $G(\PP^n,d)$), corresponding to splitting types of curves. The images in   $\ol{M}^k$ induce a canonical stratification. We index the strata by nested subsets $I$ of the power set of $D:=\{1,..., d\}$. Elements of $I$ label marked points or subcurves of the generic curve parametrized by the stratum. The incidence relations among elements of $I$ reflect incidence relations among the assigned marked points and subcurves.

\begin{notation}
Let $\cP$ denote the power set of $D$, and let $$\cP_{k}= \{ h\in \cP;  |h|=k\}$$ and similarly let $\cP_{>k}$ be the set of cardinal $> k$ subsets of $D$, $\cP_{\leq k}$, the set of cardinal $\leq k$ subsets. Similarly, for $I\subset \cP$, set $I_k:=I\bigcap \cP_k$ etc.

\end{notation}

\begin{definition}
$I\subset \cP\setminus\{\emptyset, D\}$ is a nested set if, for any two $h, h'\in I $, the intersection $h\cap h'$ is either $h$, $h'$ or $\emptyset$.
\end{definition}

\begin{definition}
Fix a positive number $k<d$ and a nested set $I\subset \cP\setminus\{\emptyset, D\}$ such that $h\cap h'=\emptyset$ for any distinct  $h, h'\in  I_{\leq d-k}$. An $I$--type, $k$--stable, degree $d$ map from a rational curve into $\PP^n$ is a tuple
 $$(C, p_1, \{p_h\}_{h\in I_{\leq d-k}}, \{C_{h}\}_{h\in I_{> d-k}}, \cL, e)$$ of a   $k$--stable, degree $d$ pointed map $(C, p_1, \cL, e)$, together with  marked points $\{p_h\}_{h\in I_{\leq d-k}}$ and connected subcurves $ \{C_{h}\}_{h\in I_{> d-k}}$ satisfying the following properties:
\begin{enumerate}
\item $\forall h\in  I_{> d-k}$,  $p_1 \not\in C_h\subset C$ and $\deg \cL_{|_{C_h}}=|h|$;
\item $\forall h\in  I_{\leq d-k}$, $\dim \Coker e_{p_h}=|h|$;
\item compatibility of incidence relations:
              \begin{itemize}    \item $\forall h\in  I_{\leq d-k}$, $\forall h'\in  I_{> d-k}$, $h\subset h'$ iff $p_h\in C_{h'}$;
\item  $\forall h, h'\in  I_{> d-k}$, if $h'\subset h$, then $C_{h'}\subset
C_{h}$, if $h\subset h'$, then $C_{h}\subset C_{h'}$, otherwise $C_{h}\cap
C_{h'}=\emptyset $; \end{itemize}  \end{enumerate}

A curve $C$ which admits a set of points and components with the above
properties is said to be of  $I$--splitting type.  \end{definition}

Intuitively, the points in $h$ may be thought of as the pullback on $C_h$ of a
hyperplane divisor in $\PP^n$.

\begin{notation}

Given any nested set $I\subset\cP$, let $G_I\subset S_d$ be the largest subgroup that keeps each $h\in I$ fixed. $G_I$ decomposes into a direct sum of permutation groups $S_{h'}$, where $h'\in \cP$ is one of the sets $h\setminus (\bigcap_{h''\in I, h''\not=h}h'')$ for $h\in I$, or $D\setminus (\bigcap_{h''\in I}h'')$.

\end{notation}

\begin{proposition} a) The functor of isomorphism families of type $I$, $k$--stable, degree
$d$ maps into $\PP^n$ is coarsely represented by a scheme $\ol{M}_{I}^k$, having a canonically associated smooth Deligne--Mumford stack $\ol{M}_{I}^k$. (This stack actually finely represents the functor).

b) Given two nested sets $I\subset J$ as above, then  $
\bar{\phi}^I_J:\ol{M}_{J}^k  \to \ol{M}_{I}^k$ is a regular local embedding of
stacks, as long as $(J\setminus I)_{< d-k}=\emptyset $. In particular,
$\ol{M}_{I}^k$ admits a regular local embedding into $\ol{M}^k:=\ol{M}_{0,1}(\PP^n,d,k)$. The image is a closed substratum $\ol{M}^k(I)$ of
$\ol{M}_{0,1}(\PP^n,d,k)$  parametrizing $k$--stable, degree $d$--maps
from curves of $I$--splitting type into $\PP^n$. 

c) If $|I\cap \cP_l|=m>1$ for some positive $l$, then the group of symmetries
$S_m$ acts nontrivially on the stack $\ol{M}_{I}^k$, permuting the elements of
$I_l$, and for any $h, h'\in I_l$, the transposition $\tau$ of $\{h, h'\}$
induces a commutative diagram $$\diagram \ol{M}_{I}^k \rto^{\tau}\dto &
\ol{M}_{I}^k \dto \\ \ol{M}_{I\setminus\{h\}}^k  \rto^{=}
&\ol{M}_{I\setminus\{h'\}}^k  \enddiagram $$

\end{proposition}

\begin{proof}

a) A finite presentation for  $\ol{M}_{I}^k$ is constructed via rigid $k$--stable maps. We fix a homogeneous coordinate system $\bar{t}$ on $\PP^n$. Recall the sequence of morphisms:
\bea \ol{M}_{0,1}(\PP^n,d,k,\bar{t}) \hookrightarrow  \begin{CD}
G(\PP^n,d,k,\bar{t}) @>{f_{\bar{t}}}>> \PP^n_d(\bar{t}) @>{p_{\bar{t}}}>>
(\PP^1)^{d(n+1)}, \end{CD} \eea  where $f_{\bar{t}}$  is the composition of blow--down
morphisms described in section 1 and  $\PP^n_d(\bar{t})$ is a
$(\CC^*)^n$--torus over an open subset of $(\PP^1)^{d(n+1)}$.

 Let $N=\{0,...,n\}$ and let $\Delta_{N\times h}$ denote the 
diagonal in $(\PP^1)^{d(n+1)}$ where the coordinates corresponding to the set $N\times h \subset
N\times \{1,...,d\}$ agree. $\Delta_{N\times h}^k(\bar{t}) \subset 
G(\PP^n,d,k,\bar{t})$ is defined inductively as follows:
$\Delta_{N\times h}^0(\bar{t}) =p_{\bar{t}}^{-1}(\Delta_{N\times h})$, and
$\Delta_{N\times h}^k(\bar{t}) $ is the strict transform of $\Delta_{N\times
h}^{k-1}(\bar{t}) $ by the $k$--th blow--up, except at the $|h|$--th
step when the total transform is considered instead. Define $D_{N\times h}^k(\bar{t}) :=
\Delta_{N\times h}^k(\bar{t}) \cap \ol{M}_{0,1}(\PP^n,d,k,\bar{t})$ and  $ D_{N\times I}^k(\bar{t}):=\bigcap_{h\in I}
D_{N\times h}^k(\bar{t})$. 
 $D^k_{N\times I}(\bar{t})$ parametrizes $I$--type, $k$--stable, 
degree $d$, $\bar{t}$--rigid maps, in the sense of Definitions 1.9
and  2.2. We will also denote this space by $\ol{M}_{I}^k(\bar{t}).$
By the
standard method employed in Proposition 4 of \cite{fulton1}, the GIT quotients
$$\left( \bigsqcup_{[g]\in (S_d/G_I)^{n+1}} D_{g(N\times
I)}^k(\bar{t})\right) / (S_d)^{n+1}\cong  D^k_{N\times I}(\bar{t})/(G_I)^{n+1}$$ 
 glue together to form the coarse moduli space $\ol{M}^k_I$.   The \'etale presentation of the
stack $\ol{M}_{I}^k$ is formed, as in Proposition 1.11, from GIT quotients of 
small neighborhoods of points $x \in \ol{M}_{I}^k(\bar{t})$ by the largest small
subgroup in the stabilizer of $x$ in $(G_I)^{n+1}$, and for various $\bar{t}$.

% de explicat prima oara cum eliminam elementele grupului care fixeaza puncte in spatiul de moduli fara a fi induse de automorphisme ale curbei reprezentate de punct.  tot acolo- cum apar subgrupurile mici: pe curbe cu puncte marcate obtinute identificand mai multe sectiuni.

 b) If $I$, $J$ are nested sets as above, and if $I\subset J$, 
then there is an obvious morphism $ \bar{\phi}^I_J: \ol{M}^k_J \to 
\ol{M}^k_I$, induced by inclusion on the $\bar{t}$--covers, as $G_J\subset G_I$. Moreover, if
$(J\setminus I)_{< d-k}=\emptyset$, then for any $x \in
\ol{M}_{J}^k(\bar{t})\hookrightarrow \ol{M}_{I}^k(\bar{t})$,  the largest small subgroup of its stabilizer in  $(G_J)^{n+1}$ coincides with that in $(G_I)^{n+1}$. By this feature the morphism $ \bar{\phi}^I_J$ is a
regular local embedding, as seen on the \'etale covers of the stacks. This is
not the case if there is an element $h\in (J\setminus I)_{< d-k}$: then
 for most $h''\in\cP_{\leq d-k}$ such that $h\subset h''$, and for $x\in D_{N\times h''}^k(\bar{t}) \cap D_{N\times J}^k(\bar{t})$, the small group $(S_{h''})^{n+1}$ is a summand in the stabilizer of $x$ in $(G_I)^{n+1}$, while 
 $ G_J\cap S_{h''}= S_{h}\oplus S_{h''\setminus h}$. This
prevents the morphism  $ \bar{\phi}^I_J$ from being a regular local embedding in
this case.

c) follows from the local structure described in a). 
\end{proof}

\begin{notation}

 Let $I\subset\cP$ be a nested set.  Then the largest finite covers of $\ol{M}^k_{I}$ which embedd into $\ol{M}_{0,1}(\PP^n,d,k,\bar{t})$ are $\bigsqcup_{[g'] \in (S_d/G_I)^n}D^{k}_{(g_0, g')(N\times I)}(\bar{t})$, for $g_0\in S_d/G_I$. This follows from the stability condition (2) in Definition 1.9. For a point $x$ of this cover and $t:=(x, \bar{t})$, let $\ol{M}^k_{I}(t)$ denote the \'etale local cover of $\ol{M}^k_I$ obtained as a quotient by the appropriate small group of a neighborhood of $x$ in the above finite cover, as in Proposition 1.11. The \'etale covers $\{\ol{M}^k_{I}(t)\}_t$ organize the stacks $\{\ol{M}^k_I\}_I$ into an $S_d$-- network of regular local embeddings, which shall be discussed in detail in section 3.

\end{notation}

 Consider a nested subset $I\subset\cP$. At the level of local $(G_I)^{n+1}$-- covers,
  the morphism $f^k_I(\bar{t}): \ol{M}_I^{k+1}(\bar{t})\to\ol{M}_I^k(\bar{t})$ factors 
into a sequence of regular blow--ups and projections of projective bundles. The blow--ups 
are along $\ol{M}_{Ih}^k(\bar{t})$, for all $h\in \cP_{d-k}\setminus I$. The projective
 bundles are restrictions of the exceptional divisors in the blow--ups of 
$\ol{M}^k(\bar{t})$ along $\ol{M}^k_{h}(\bar{t})$, for all $h\in I_{d-k}$. Their 
quotients make up projective fibrations over $\ol{M}_I^k$:
\begin{notation}

 Let $\PP_h:=\ol{M}_{0,1}(\PP^n,|h|,1)\times_{\PP^n}\ol{M}^k_h$, where the fiber product 
is taken along evaluation functions. Note that there is a natural evaluation function at the $h$--th point $ev_h:\ol{M}^k_h\to \PP^n$.
 Indeed, for a generic point $(C, p_1, p_h, \cL, e)$ of $\ol{M}^k_h$, $e$ induces a linear
 system on $\cL(|h|p_h)$ whose base locus never contains $p_h$. There is also a natural 
evaluation function $ev_1:\ol{M}_{0,1}(\PP^n,|h|,1)\to\PP^n$. 

 $\PP_h$ is obtained by gluing $(G_h)^{n+1}$--quotients of the projective bundles
 $\PP(\cN_{\ol{M}^k_{h}(\bar{t})|\ol{M}^k(\bar{t})})$. It will play an important role 
in the later Chow ring computations. 

 For any $J\subset\cP_{d-k}$, let $\PP_J$ denote the fiber product of all pullbacks to $\ol{M}_J^k$  of the fibrations $\PP_h$. The morphism $f^{k}_{J}: \ol{M}_{J}^{k+1}\to \ol{M}_J^k$ factors into $g_J:  \ol{M}_{J}^{k+1}\to\PP_J$ and $\pi_J: \PP_J\to \ol{M}_J^k$. 

\end{notation}

 Reversing the order, one could also think of $\ol{M}_J^{k+1}$ as a projective fibration. There is a fibre square:
\bea \diagram \ol{M}_J^{k+1}\dto \rto & \PP_J \dto   \\
                \ol{M}^{k+1}(^k_J) \rto & \ol{M}_J^{k} \enddiagram \eea

Here $\ol{M}^{k+1}(^k_J)$ is the space of $(k+1)$--stable, degree $(d-|J|(d-k))$--maps 
from rational curves into $\PP^n$, with a set of distinct, smooth marked points 
$\{p_1, \{p_h\}_{h\in J} \}$. Locally this  is the $(G_J)^{n+1}$-- quotient of the blow--up 
of $\ol{M}_{J}^k(\bar{t})$ along all $ \ol{M}^{k}_{Jh}(\bar{t})$, for $h\in \cP_{d-k}$ 
disjoint from the elements of $J$.

\begin{remark}
Consider any $h\in\cP_{d-k}\setminus \{\emptyset, D\}$ and any nested set $I$ such
that $I\not= I\cup \{h\}$ is again a nested set. The space $\ol{M}^k_{ Ih}$ maps 
into one component of the blow--up locus at the $k$--th step in $\ol{M}^{k}_I$. The components of the blow--up locus correspond to the different types of positionings of 
cardinal $k$ subsets $h\subset D$ with respect to the elements of $I$.                               
                                                                             
\end{remark}

%%partitia corespunzatoare unui cuib de multimi
%spatiul de moduli= produs de spatii fibrate.

\section{The Chow ring }

  Previously we have found $\ol{M}_{0,1}(\PP^n, d)=:\ol{M}^d$ to
be the end product of a sequence of birational transformations. This section is concerned with the Chow ring computation of this
space.  The inductive method centers on
writing the Chow ring of  $\ol{M}^{k+1}_I$ as an algebra over that of
$\ol{M}^{k}_I$, for any nested subset $I\subset \cP$ and any step $k$. The
calculation is done in two parts, each adapting methods employed for a usual
blow--up (see \cite{fulton3}, \cite{keel}): the classical calculation by
Grothendieck of the Chow ring for a projective fibre bundle is extended to
the case of a weighted projective fibration via localization. Then, the Chow
ring of a weighted blow--up along a regular local embedding is
 related to that of its exceptional divisor via essentially a self--intersection
formula. As the exceptional divisor displays multiple components at the level
of \'etale covers, its self intersection becomes more involved and justifies
the introduction of a network containing multiple copies of strata. 

\subsection*{Localization}
Here the Chow ring of a weighted projective fibration is computed as
an algebra over the Chow ring of the basis, via localization. The
working examples are $\ol{M}^1_I\to \PP^n$ and 
$\PP_h\to \ol{M}^k_h$. In the first case, the underlying action of the group $(G_I)^{n+1}$
 on a cover plays an important role, in the second case this group is $(S_h)^{n+1}$. 
 The first lemma defines the properties of a weighted blow--up, going from local to global presentation.

  Let $G$ be a finite group.

\begin{lemma}
  Consider a  birational morphism of reduced schemes
 $\pi : \tilde{X}\to X$. Suppose that $X$ can be covered by open sets of type $U\cong
 U'/G$  and that $\pi^{-1}(U)\cong \tilde{U}'/G$, where $ \tilde{U}'$
 is the blow--up of $U'$ along a reduced locus $Y'$ which is left point--wise invariant by the action of $G$. Then 

a) $\tilde{X}$ is a weighted blow--up of $X$
 along a locus $Y$, namely, $$\tilde{X}=\Proj (\oplus_{n\geq 0} \cI_n)$$ for an increasing filtration $\{  \cI_n\}_{n\geq 0}$ of the ideal $\cI_Y$, such that  $\cI_0=\cO_X$,  $\cI_1=\cI_Y$ and $\cI_n\cI_m\subseteq\cI_{m+n}$ for all $m,n\geq 0$.

b) The reduced structure of the exceptional divisor in $\tilde{X}$ is $$E=\Proj(\oplus_{n\geq 0} \cI_n/\cI_{n+1}). $$ 

 \end{lemma}

\begin{proof}

a) Let $E$ denote the exceptional divisor of $\pi: \tilde{X} \to
X$, and $Y$ its image in $X$.  $\cI_E$ and $\cI_Y$ will denote their sheaves of ideals. 
The sheaves $\cI_n:= \pi_*\cI_E^n$ form a
filtration of $\cI_Y$ such that $\cI_0=\cO_X$
and $\cI_j\cI_k\subseteq \cI_{j+k}$ for all indices $j,k$.  It remains to show that
$\tilde{X}=\Proj (\oplus_{n\geq 0} \cI_n)$. Indeed,
$\pi^{-1}(U)=\Proj(\oplus_{n\geq 0}(f_*\cI_{Y'}^n)^G)$, where $f$
is the morphism from $U'$ to $X$. Also, if $\tilde{f}$ is the
morphism from $\tilde{U}'$ to $ \tilde{X}$, then over
$\pi^{-1}(U)$, $\cI_E^n\cong (\tilde{f}_*\cI_{E'}^n)^G,$ where $E'$ is the exceptional divisor of $\tilde{U}'$. If
$\pi': \tilde{U}'\to U'$ is the blow--down morphism, then
 $(f_*\cI_{Y'}^n)^G\cong (f_*\pi'_*\cI_{E'}^n)^G\cong \pi_*(\cI_E^n)_{|U}\cong \cI_{n_{|U}}$.

In the same way,  $\cI_{n}/\cI_{n+1}\cong \pi_*(\cI_E^n/\cI_E^{n+1})$.

b) We claim that the reduced algebra of $\cO_Y\otimes (\oplus_{n\geq 0} \cI_n)$
is $\oplus_{n\geq 0} \cI_n/\cI_{n+1}$.
 Let $\cJ$ be the ideal of $\cO_Y\otimes (\oplus_{n\geq 0} \cI_n)$
generated by the images of the morphisms $\phi_k:\cO_Y\otimes
\cI_{k+1}\to\cO_Y\otimes \cI_{k}$. The ideal $\cJ$ is nilpotent.
Indeed, let $x$ be a section of Im $\phi_k$. For any $j>k$,
$x^j=x^{j-1}x$ is a section of $\cI_{(k+1)(j-1)}\cI_Y\subset
\cI_{kj}\cI_Y$. On the other hand, $\cO_Y\otimes (\oplus_{n\geq 0}
\cI_n)/\cJ\cong \oplus_{n\geq 0} \cI_n/\cI_{n+1}$ and
$\oplus_{n\geq 0} \cI_n/\cI_{n+1}$ is reduced because it is
locally the $G$--invariant part of $\oplus_{n\geq 0}
\cI_{Y'}^n/\cI_{Y'}^{n+1}$.

\end{proof}

 Keeping the notations from above, further assume that each 
  $U'$ is a smooth scheme on which $G$ acts as a small
 group, and that $Y'\subset U'$ is a regular embedding.
 This local situation may be expressed globally by the following properties of the filtration $\cI_k$:
\begin{enumerate}
\item  $\cI_k\cap \cI_Y^2=\sum_{j=1}^{k-1} \cI_j\cI_{k-j}$,
\item $ \cI_k/(\cI_k\cap\cI_Y^2)$ is a subbundle of the conormal bundle $\cI_Y/\cI_Y^2$.
\end{enumerate}

These requirements are the minimum that insure  a natural structure of smooth Deligne--Mumford stack on $\tilde{X}$. Under these conditions, the following Lemma holds:

\begin{lemma}

a) The normal bundle in $A:=\Spec(\oplus_{n\geq 0}\cI_n/\cI_{n+1})$ of the fixed locus $Y$ under the natural $\CC^*$ action on $A$ is
$$\cN_{Y|A}= \oplus_{n\geq 1}\cN_n/\cN_{n+1},$$
where $\{\cN_n\}_n$ is the filtration of the normal bundle $\cN_{Y|X}$ dual to the filtration $ \{\cI_n/(\cI_n\cap\cI_Y^2)\}_n $ of $\cI_Y/\cI_Y^2$.

b) There is a ring isomorphism $$A^*(E; \QQ)\cong \frac{A^*(Y; \QQ)[\tau ]}{ < P_{Y|X}(\tau )>},$$
where $P_{Y|X}(t)$ is the top equivariant Chern class of the bundle $ \cN_{Y|A}$. In particular, the free term of $P_{Y|X}(t)$ is the top Chern class of $\cN_{Y|X}$. $\tau$ is the first Chern class of $\cO_E(1)$. 

\end{lemma}
\begin{proof}

a)The ideal of the zero section $Y$ in $A$ is $\oplus_{n\geq 1}
\cI_n/\cI_{n+1}$. There is a natural morphism of modules over
$\oplus_{n\geq 0}\cI_n/\cI_{n+1}$:
 $$\oplus_{n\geq 1}\cI_n/\cI_{n+1}\to \oplus_{n\geq 1}\cI_n/(\cI_{n+1}+\cI_n\cap
 \cI_Y^2)\cong \oplus_{n\geq 1}\frac{\cI_n/\cI_{n}\cap
\cI^2_Y}{\cI_{n+1}/\cI_{n+1}\cap \cI^2_Y}.$$ The kernel is
$\oplus_{n\geq 1}(\cI_n\cap\cI^2_Y)/(\cI_{n+1}\cap\cI^2_Y)$, which
by property (1) above, is isomorphic to $(\oplus_{n\geq
1}\cI_n/\cI_{n+1}))^2$.

b) Consider a finite $G$-- cover $\{ U'_i\}_i$ of $X$ and $\{ Y'_i\}_i$ of $Y$ such that 
the blow--ups of $U'_i$ along $Y'_i$ form a $G$--cover of $\tilde{X}$. Let $\{E_i\}_i$ 
 form the corresponding open cover of $E$. A canonical stratification $\{V_i\}_{ i}$ of $E$ is assigned: $V_i=Z_{i-1}\backslash Z_i$ where
$Z_0:=E$ and $Z_i:=E\setminus (\bigcup_{j=1}^{i}E_{j})$ for
$i\geq 1$. A standard argument based on the open--closed exact sequence of Chow
groups of each $Z_{i-1}=V_i\cup Z_{i}$ leads to the following generators of
$A^*(E)$ as an $A^*(Y)$--module: $\tau^i$, for all $i$ such that $0\leq i < e:=\codim_XY$. Here $\tau$ denotes $c_1(\cO_{E}(1))$.

 Moreover the following relations hold: $\pi_{
*}(\tau^i)=0$ for all $i<e$ and $\pi_{*}(\tau^e)=[Y]$, where $\pi$ is the morphism
$\pi : E\to Y$. We will now find a relation 
\bean \tau^e+\sum_{i=0}^{e-1} a_i \tau^i=0 \eean on
$E$. There is a $\CC^*$--equivariant morphism
\bea F: N_{E|\tilde{X}}:=\Spec (\oplus_j \cI_E^n/\cI_E^{n+1}) \to A=\Spec(\oplus_{n\geq 0}\cI_n/\cI_{n+1}) \eea
for the natural $\CC^*$ -- actions having $E$ and $Y$ as fixed loci. 
An elementary form of Atiyah--Bott localization formula applied to $F$ yields: \bea
F_{*}  \frac{1}{ c^{eq}_{top}(  \cN_{E|\tilde{X}} )}  =
\frac{1}{ c^{eq}_{top} ( \cN_{Y|A} )}, \eea
where $c^{eq}_{top}$ denotes the equivariant top Chern class.
 As $c^{eq}_{top}(\cN_{E|\tilde{X}})= t-\tau$, it follows that
$c^{eq}_{top}(\cN_{Y|A})$ is the inverse of the Segre series
$\sum_{j\geq e-1} F_{*}(\tau^{j}) t^{-j-1}.$ This implies that $\tau$ is
the root of the polynomial $c^{eq}_{top}(\cN_{Y|A})$. 
\end{proof}

 The above lemmas apply to any projective morphism of stacks $P\to Y$ together with a  sheaf $\cL$ on $P$, satisfying the following conditions: $Y$ may be covered by open subsets $Y'$ such that there exist a $G$--equivariant projective
 bundle $\pi': \PP(V')\to Y'$ and a morphism $\tilde{g} : \PP(V')\to P$
 whose image is $\PP(V')/G$; locally $\cL$ is the sheaf of $G$--invariant sections of  $\cO_{\PP (V')}(1)$ and $G$ acts as a
 small group on $ \pi'_*\cO_{\PP (V')}(1) $. Then $P\cong \Proj (\oplus_n \pi_* \cL^n)$ is the exceptional divisor of the following weighted blow--up:
$$\tilde{X}:=\Spec(\oplus_n\cL^n)\to X:= \Spec (\oplus_n \pi_* \cL^n).$$
 We say that $P$ is a weighted projective fibration.

The first application of interest to us is to the first intermediate moduli space and 
its substrata $\ol{M}_I^{1}$, for a partial partition $I$:  
\begin{notation}

  A set $I=\{h_1,...,h_s\}$ is called a partial partition of
$\{1,...,d\}$ if ${h_1},...{h_s} \subset \{1,...,d\}$ are disjoint subsets.    
We denote by  $l_I:=\sum_{i=1}^s |{h_i}|.$

\end{notation}

There are two important classes on $\ol{M}_I^{1}$: the pullback $H$ of the hyperplane 
class in $\PP^n$ via the natural evaluation map; and the cotangent class $\psi$, the first
 Chern class of the line bundle $s_1^*(\omega_{\ol{U}_I^{1}|\ol{M}_I^{1}})$, where 
$\pi:\ol{U}_I^{1}\to \ol{M}_I^{1}$ is the universal family and $s_1$ is its canonical 
section.

\begin{lemma}
The Chow ring with rational coefficients of $\ol{M}_I^{1}$ is
$$ A^*(\ol{M}_{I}^{1};\QQ)=\frac{\QQ[H,\psi]}{\langle
H^{n+1}, P_{\ol{M}_{I}^{0}| \ol{M}^0}(\psi) \rangle} $$
 where
\bea
P_{\ol{M}_{I}^{0}|\ol{M}^{0}}(t)=t^{s-1}\prod_{i=1}^{{d-l_I}}(H+it)^{n+1}. \eea
\end{lemma}

\begin{proof}
Recall that
$ G(\PP^n,d,1)$ is the weighted blow--up of $\PP^n_d$ along
$\PP^1\times\PP^n$ in the sense of Lemma 3.1, and
 $\ol{M}^1$ is the restriction of the exceptional
divisor to $\{ 0\}\times\PP^n$. Given any partial partition $I$, $\ol{M}^1_I$
is the restriction to $\{ 0\}\times\PP^n$ of 
exceptional divisor $E$ in the weighted blow--up of
$(\PP^1)^{s}\times\PP^n_{d-l_I}$ along $\PP^1\times\PP^n$.
By a classical argument, pullback of the class $[E]$ from this weighted blow--up to $\ol{M}^1_I$ is $-\psi$.
 If the points of $\PP^n_d, \PP^n_{d-l_I}$ are understood 
modulo constants as $(n+1)$--tuples of polynomials in one variable $t$, then
the map $$(\PP^1)^{s}\times\PP^n_{d-l_I}\to\PP^n_d$$  is the multiplication of
each of the $(n+1)$ degree $(d-l_I)$ polynomials by the polynomial
$\prod_{i=1}^s(b_it-a_i)^{|h_i|}$, for each $\{(a_i:b_i)\}_{1\leq i\leq
s}\in (\PP^1)^s$. The blow--up locus $\PP^1\times\PP^n$ is embedded by:
 $$\PP^1\times\PP^n\to (\PP^1)^s\times \PP^n_{d-l_I}$$
 $$((a:b), (g_0:...:g_n))\to
 (\{(a:b)\}_{1\leq i\leq s},((at-b)^{d-l_I}g_0:...:(at-b)^{d-l_I}g_n))). $$
The normal bundle $$\left.\cN_{\PP^1\times\PP^n |
(\PP^1)^{s}\times\PP^n_{d-l_I}}\right|_{\PP^n}=
(\cO_{\PP^n}^{s}\oplus\oplus_{i=1}^{(n+1)(d-l_I)}\cO_{\PP^n}(1))/\cO_{\PP^n}$$
admits a natural filtration  $\{\cN_k\}_k$ described in Lemma 3.2. Here
$$\cN_k=(\cO_{\PP^n}^{s}\oplus\oplus_{i=1}^{(n+1)(d-l_I-k)}\cO_{\PP^n}(1))/\cO_{\PP^n}=\left.\cN_{\PP^1\times\PP^n|
(\PP^1)^{s}\times\PP^n_{k}}\right|_{\PP^n}.$$
   $\CC^*$ acts on the bundle
$\cN=\oplus_k \cN_k/\cN_{k+1}$  with weights $(1,...,d-l_I)$, and the
top equivariant Chern class of $\cN$ is thus
$$P_{\PP^1\times\PP^n|(\PP^1)^{s}\times\PP^n_{d-l_I}}(t)=t^{s-1}\prod_{i=1}^{d-l_I}(H+it)^{n+1}. $$
Thus by Lemma 3.2,
 \bea
A^*(\ol{M}_{I}^{1};\QQ)=\frac{\QQ[H,\psi]}{\langle H^{n+1},
P_{\PP^1\times\PP^n|(\PP^1)^{s}\times\PP^n_{d-l_I}}(\psi) \rangle}. \eea
 
 Looking at the pull--back morphism $A^*(\ol{M}^1;\QQ)\to
A^*(\ol{M}_{I}^{1};\QQ)$ we find the class of $\ol{M}_{I}^{1}$
in $\ol{M}^1$ to be a multiple of
$\psi^{-s}\prod_{d-l_I+1}^{{d}}(H+i\psi)^{n+1}$. For dimension
reasons, the coefficient is a constant and moreover, after
push--forward to $\PP^n$ we find the constant to be 1. Therefore
\bean
[\ol{M}_{I}^{1}]=\psi^{-s}\prod_{d-l_I+1}^{{d}}(H+i\psi)^{n+1}.\eean

\end{proof}

The next example where Lemma 3.2 applies is the fibration $\PP_n\to \ol{M}^k_h$ defined in section 2. Indeed, $\PP_h$ is the pullback to $\ol{M}^k_h$ of the weighted projective fibration $\ol{M}_{0,1}(\PP^n,|h|,1)\to\PP^n$, and therefore, its Chow ring is generated over $A^*(\ol{M}^k_h)$ by the class $\psi_h$, the pullback of the cotangent class on $\ol{M}_{0,1}(\PP^n,|h|,1)$.
 From a different perspective, the conormal bundle of $\ol{M}^{k+1}_h$ in $\ol{M}^{k+1}$ 
descends to $\PP^n$, being constant on the fibres of $g_h:\ol{M}^{k+1}_h  \to \PP_h$. The first Chern  class $\tau_h$ of the descent bundle is an alternate generator of the Chow ring. Moreover,

\begin{lemma} There is a filtration $\{\cN_j\}_j$ of the normal bundle $\cN_{\ol{M}^k_h|\ol{M}^k}$ such that
$$A^*(\PP_h; \QQ)= A^*(\ol{M}^k_h)[\tau_h]/ < P_h(\tau_h) >, $$
where $P_{h}(t)$ is the top equivariant Chern class of the fibre bundle $\cN_h=\oplus_j \cN_j/\cN_{j+1}$ with respect to the natural $\CC^*$ action.
\end{lemma}
 
\begin{proof}
 Consider the $\bar{t}$--covers of sections 1 and 2, and the blow--up $M(\bar{t})$ of $\ol{M}^k(\bar{t})$ along $\ol{M}^k_h(\bar{t})$. The quotient by $(G_h)^{n+1}$ of the exceptional divisor  is an open subset of $\PP_h$. Following the discussion in Proposition 1.11, quotients by small subgroups in $(G_h)^{n+1}$ yield locally an \'etale cover of $\PP_h$. 

The quotient $M(\bar{t})/(S_h)^{n+1}$ is a weighted blow--up and thus Lemma 3.2 applies 
to the \'etale covers of $\PP_h$, yielding the  filtration $\{\cN_j\}_j$ of the normal 
bundle $\cN_{\ol{M}^k_h|\ol{M}^k}$. Indeed, for all $x\in M(\bar{t})$, the largest small 
subgroup $H_x$ of the stabilizer of $x$ decomposes into $(S_h)^{n+1}\oplus H'_x$, such that the weighted blow--up structure of $M(\bar{t})/(S_h)^{n+1}$ descends through $H'_x$. Functoriality of the constructions in Lemmas 3.1 and 3.2 insures the gluing property of the filtration $\{\cN_j\}_j$ under different \'etale maps, and for different $\bar{t}$--s.
\end{proof}

The preceding examples set up the first induction steps in the calculation of $A^*(\ol{M}^k_I)$. The general step relies on the following construction:

\subsection*{The Chow ring associated to a network of regular local embeddings}

Let $\cP$ be a finite set and let $G$ be a finite group acting on $\cP$. 
% cum actioneaza G pe $\cP$?

For every subset $I\subset \cP$,  $G_I$ denotes the largest
 subgroup of $G$ that fixes all elements of $I$. Among subsets of $\cP$, we consider a special family closed under the $G$--action, which will be called the family of      allowable sets, such that if  $I$ is allowable, then any of its subsets is allowable too.

\begin{definition}
 A $G$-- network of regular local embeddings indexed by allowable subsets $I$ of $\cP$
is a category of smooth Deligne--Mumford stacks $\{\ol{M}_I\}_{I\subset \cP }$
with unique morphisms $\bar{\phi}^I_J :\ol{M}_J\to \ol{M}_I$ for every
$J\supset I$, canonical isomorphisms $g: \ol{M}_I\to \ol{M}_{g (I)}$
for every $[ g ]\in G/G_I$, and a set of \'etale covers by
schemes $\ol{M}_I(t)\to \ol{M}_I$, such that
there is a Cartesian diagram
\bea  \diagram  \bigsqcup_{[ g ]\in G_I/G_J} \ol{M}_{g(J)}(t)
\rto \dto & \ol{M}_I(t)\dto \\
\ol{M}_J \rto^{\bar{\phi}^I_J } & \ol{M}_I,
\enddiagram \eea
each morphism $\ol{M}_J(t)\to \ol{M}_I(t)$ is an inclusion, $\ol{M}_J(t)\bigcap \ol{M}_K(t)=\ol{M}_{J\bigcup K}(t)$ in $\ol{M}_{J\bigcap K}(t)$, and all intersections are transverse. \end{definition}

We will write $IJ$ for $I\bigcup J$ when $I\cup J$ is allowable.

The stacks $\ol{M}_I$ with $I\not=\emptyset$ minimal among the allowable
subsets of $\cP$ will be called generators of the network. 

\begin{notation}
 We will denote the image of $\bar{\phi}^I_J$ by $\ol{M}_I(J\setminus I)$.
This depends only on the class of $J$ modulo the action of
elements in $G/G_I$. The morphism $\bar{\phi}^I_J$ factors through a finite map
$\phi^I_J:\ol{M}_J\to  \ol{M}_I(J\setminus I)$ and an embedding $j^I_J: \ol{M}_I(J\setminus
I)\to \ol{M}_I.$ When $I=\emptyset$ we omit the subscript $I$. Thus $\ol{M}$
denotes the final element of the category.  \end{notation}

 A notion of extended Chow ring is assigned to a $G$-- network of regular local
embeddings by concatenating the Chow rings of all network elements, modulo a
natural equivalence relation which keeps track of the structure at the level
of \'etale atlases: 

\begin{notation} Fix $I\subset \cP$ allowable and $h\in \cP$ such that $I\cup \{h\}$ is allowable. 
 To any cycle $\alpha=[V] \in Z_l(\ol{M}_I(h))$  we will associate a
cycle $\alpha_h \in  Z_l(\ol{M}_{Ih})$ as follows: $$\alpha_h=\frac{\sum_i
 [V^i_h]}{deg\left( (\phi^{I}_{Ih})^{-1}(V)/V\right) }$$
where $\{V^i_h\}_i$ are the $l$--dimensional components of $(\phi^{I}_{Ih})^{-1}(V)$.
\end{notation}
\begin{definition}
 The graded $\QQ$--vector space $B^*(\ol{M};\QQ)$ is defined as follows:$$B^*(\ol{M};\QQ):= \bigoplus_{l=0}^{\dim
(\ol{M})}B^l(\ol{M}),$$ where the extended Chow groups are
\bea   B^l(\ol{M})  :=\oplus_{I} Z^{l-\codim_{\ol{M}}\ol{M}_{I}}(\ol{M}_{I}) / \sim ,  
\eea  the  sum taken after all allowable subsets $I\subset \cP$ with $
\codim_{\ol{M}}\ol{M}_I \leq l$. The equivalence relation $\sim $ is
generated by rational equivalence together with relations of the type: $$
j^I_{Ih *}\alpha \sim \sum_{[g ]\in G/G_h}\alpha_{g (h)}, $$ for any
cycle $\alpha\in Z_l(\ol{M}_I(h))$, where $\alpha_{g (h)}$ is the
corresponding cycle in $\ol{M}_{Ig (h)}$. \end{definition}

\begin{example} With the notations from section 2, 
 $G:=S_d$ acts on $\cP_{> d-k}:=\{h\subset\{1,...,d\};  |h|> d-k\}$. Whenever
$l\geq k-1$, the set of strata $\ol{M}^l_I$ indexed by nested
subsets $I\subset \cP_{> d-k}$ forms an $S_d$-- network. The associated
extended Chow ring will be denoted by $ B^*(\ol{M}^l, k)$. For example, $
B^*(\ol{M}^k, k)$ corresponds to the network generated by normalizations of all exceptional
divisors (and their strict transforms) from the steps 1 to $k$. 

Pullbacks by the
morphisms $f^l_I:\ol{M}_I^{l+1}\to \ol{M}_I^{l}$ are compatible with the
equivalence relation $\sim$ and add up to a global pullback $f^{l *}: 
B^*(\ol{M}^l, k)\to  B^*(\ol{M}^{l+1}, k)$.

  \end{example}

\begin{example}
 Given any $G$--network and a fixed member $\ol{M}_I$, the sets
$\{\ol{M}_J\}_{J\supseteq I}$ and their morphisms form a $G_I$--network. There
is a natural morphism of graded vector spaces:
$B^*(\ol{M}_I)\to B^*(\ol{M})$, of degree equal to the codimension of
$\ol{M}_I$ in $\ol{M}$. In addition, for $J\supset I$ there are natural
pullback and pushforward morphisms $\bar{\phi}^{I *}_J: B^*(\ol{M}_I)\to
B^*(\ol{M}_J)$ and  $\bar{\phi}^{I}_{J *}: B^*(\ol{M}_J)\to
B^*(\ol{M}_I)$ obtained  by concatenating the usual pullbacks and
pushforwards on substrata. Compatibility with $\sim$ is easily
checked. \end{example}

A distinctive feature coming from the \'etale structure of a $G$--network is 
the following ``excess intersection'' formula:
\begin{lemma}
Given any allowable  $I$, $h\in I$, and $\alpha \in A^*(\ol{M}_I(h))$, the
following relation holds: \bean \bar{\phi}^{I*}_{Ih}(j^I_{Ih *}(\alpha))=
\alpha_h\cdot
c_{top}(\cN_{\ol{M}_h|\ol{M}})+\bar{\phi}^{Ih}_{Ihh'
*}\bar{\phi}^{Ih' *}_{Ihh'}(\alpha_{h'})\eean
 where $h'=g' (h)$ for some $[g']\in G/G_h$.  Here $c_{top}(\cN_{\ol{M}_h|\ol{M}})$ denotes the top Chern class of the given normal bundle.  The term $\bar{\phi}^{Ih}_{Ihh'
*}\bar{\phi}^{Ih' *}_{Ihh'}(\alpha_{h'})$ does not depend on the
choice of $h'$. This term  supported on $\ol{M}_{Ih}(h')\subset
\ol{M}_{Ih}$ satisfies:
 \bean \bar{\phi}^{Ih}_{Ihh'
*}\bar{\phi}^{Ih'
*}_{Ihh'}(\alpha_{h'})\sim\sum_{[g]\not= [e] \in G/G_h}\bar{\phi}^{Ig(h)
*}_{Ih g(h)}(\alpha_{g(h)}). \eean
\end{lemma}
\begin{proof}
  Given an $l$--dimensional class
  $j^I_{Ih *}(\alpha)=[V]$ on $\ol{M}_I$ and the preimage $V_h$ of $V$ in  
$\ol{M}_{Ih}$, the image of $[V]$ through the  Gysin map
$\bar{\phi}^{I*}_{Ih}$ is the intersection class of  the zero section
in $\cN_{\ol{M}_{Ih}| \ol{M}_{I}}$ with   the cone of $V_h$ in $V$.
By Definition 3.5, at  the level of \'etale covers this cone has multiple
components indexed by $[g]\in G/G_h$, which yield the above sum.
\end{proof}

%acoperire etala sau embeding? de comentat.

\begin{remark}
 For any allowable $\{ h\} \cup J\subset \cP$, there is an isomorphism of normal bundles 
$\cN_{\ol{M}_{Jh}|\ol{M}_J}\cong \bar{\phi}^{h*}_{Jh}\cN_{\ol{M}_h|\ol{M}}$.
Thus for any allowable set 
$I\subset \cP$,  $$\cN_{\ol{M}_I|\ol{M}}\cong \bigoplus_{h\in
I}\bar{\phi}^{h*}_I\cN_{\ol{M}_h|\ol{M}}.$$ \end{remark}

 The next definition and lemma introduce a ring structure on
 $B^*(\ol{M})$. 

\begin{definition}
Multiplication on $B^*(\ol{M})$ is defined as follows: given any two
classes $\alpha \in A^*(\ol{M}_I)$ and $\beta \in A^*(\ol{M}_J)$, let
$$ \alpha\cdot_r\beta := \bar{\phi}^{I*}_{I\cup J}(\alpha) \cdot \bar{\phi}^{J*}_{I\cup
J}(\beta) \cdot c_{top}(\bar{\phi}^{I\cap J*}_{I\cup J}\cN_{\ol{M}_{I\cap
J}|\ol{M}})$$ in $A^*(\ol{M}_{I\cup J})$. Here $\bar{\phi}^{I*}_{I\cup J},
\bar{\phi}^{J*}_{I\cup J}$ are the (generalized) Gysin homomorphisms, as defined in \cite{vistoli}. 
\end{definition}
\begin{lemma}
 $B^*(\ol{M})$ with the multiplication $\cdot_r$ admits a
natural  structure of graded $\QQ$-- algebra.
\end{lemma}
\begin{proof}
We will denote $c_{top}(\cN_{\ol{M}_h|\ol{M}})$ by $\xi_h$.

Clearly the multiplication preserves the grading. The associativity is
 straightforward: for any  $\alpha \in
A^*(\ol{M}_I)$, $\beta \in A^*(\ol{M}_J)$ and  $\gamma \in
A^*(\ol{M}_K)$,
\bea & \alpha\cdot_r(\beta\cdot_r\gamma)=(\alpha\cdot_r\beta)\cdot_r\gamma=& \\
 &=\alpha\cdot\beta\cdot\gamma\cdot(\prod_{h\in
I\cap J\cap K}\bar{\phi}^{h*}_{I\cup J\cup K}\xi_h)\cdot(\prod_{h\in (I\cap J)\cup
(I\cap K)\cup(J\cap K)}\bar{\phi}^{h*}_{I\cup J\cup K}\xi_h).& \eea
 The compatibility of the equivalence relation $\sim$ with the
multiplication is a consequence of the previous lemma:

 Consider first $j^I_{Ih *}\alpha, \beta \in A^*(\ol{M}_I)$. Via projection
formula: \bea &j^I_{Ih *}\alpha\cdot_r\beta = j^I_{Ih *}\alpha \cdot \beta\cdot
c_{top}(\cN_{\ol{M}_I|\ol{M}})=j^I_{Ih *}(\alpha \cdot j^{I*}_{Ih}( \beta\cdot
c_{top}(\cN_{\ol{M}_I|\ol{M}}))) \sim & \\
&\sim\sum_h (\alpha \cdot j^{I*}_{Ih}(\beta\cdot
c_{top}(\cN_{\ol{M}_I|\ol{M}})))_h &
\eea
 Since $\beta\cdot c_{top}(\cN_{\ol{M}_I}|\ol{M})$ is a class on $\ol{M}_I$,
the following holds:
  $$(\alpha \cdot j^{I*}_{Ih}(\beta\cdot
c_{top}(\cN_{\ol{M}_I|\ol{M}})))_h = \alpha_h\cdot \bar{\phi}^{I*}_{Ih}(\beta\cdot
c_{top}(\cN_{\ol{M}_I|\ol{M}})) =\alpha_h \cdot_r \beta.$$

 If  $j^I_{Ih *}\alpha \in A^*(\ol{M}_I)$ and $\beta \in A^*(\ol{M}_J)$ with
$I\not=J$,  then by definition:
 \bea  & j^I_{Ih *}\alpha\cdot_r\beta = \bar{\phi}^{I*}_{I\cup J}(j^I_{Ih *}\alpha) \cdot \bar{\phi}^{J*}_{I\cup
J}(\beta) \cdot c_{top}(\bar{\phi}^{I\cap J*}_{I\cup J}\cN_{\ol{M}_{I\cap
J}|\ol{M}}).\eea
 Set $\gamma= \bar{\phi}^{J*}_{I\cup J}(\beta) \cdot
c_{top}(\bar{\phi}^{I\cap J*}_{I\cup J}\cN_{\ol{M}_{I\cap J}|\ol{M}})$. 
The morphism $\bar{\phi}^{I*}_{I\cup J}$ may be split into a composition of
pullbacks $\bar{\phi}^{I_s*}_{I_{s+1}}$, with $I_0:=I$, $I_{|J\setminus
I|}=I\cup J$ and $|I_{s+1}\setminus I_s|=1$. After successive applications of
Lemma 3.9:
\bea & \bar{\phi}^{I*}_{I\cup J}j^I_{Ih *}\alpha = \sum_{h\in
J\backslash I} \bar{\phi}^{Ih*}_{I\cup J}\alpha_h\cdot\xi_h+ \bar{\phi}^{I\cup
J}_{(I\cup J)h' *}\bar{\phi}^{Ih'*}_{(I\cup J)h'}\alpha_{h'}
\eea
 for some $h'=g(h)\not\in I\cup J$. By formula (3.4),
\bea
& \bar{\phi}^{I\cup
J}_{(I\cup J)h' *}\bar{\phi}^{Ih'*}_{(I\cup J)h'}\alpha_{h'}
\cdot \gamma= j^{I\cup J}_{(I\cup J)h'*}\left( \phi^{I\cup
J}_{(I\cup J)h' *}\bar{\phi}^{Ih'*}_{(I\cup J)h'}\alpha_{h'}
\cdot j^{I\cup J *}_{(I\cup J)h'}\gamma \right) \sim &\\
&\sim \sum_{h'=g(h)\not\in I\cup J} \bar{\phi}^{Ih'*}_{(I\cup J)h'}\alpha_{h'}
\cdot \bar{\phi}^{I\cup J *}_{(I\cup J)h'}\gamma &
 \eea
 Putting it all together, $j^I_{Ih *}\alpha\cdot_r\beta $ is equivalent to the following sum:
\bea &  \sum_{h'=g(h)\in J\backslash I}
\bar{\phi}^{Ih'*}_{I\cup J}(\alpha_{h'}) \cdot \bar{\phi}^{J*}_{I\cup J}(\beta)
\cdot c_{top}(\bar{\phi}^{(I\cap J)h'*}_{I\cup J}\cN_{\ol{M}_{(I\cap
J)h'}|\ol{M}})+&\\ & + \sum_{h'=g(h)\not\in J\cup I} \bar{\phi}^{Ih'*}_{(I\cup
J)h'}(\alpha_{h'}) \cdot \bar{\phi}^{J*}_{(I\cup J)h'}(\beta)\cdot
c_{top}(\bar{\phi}^{I\cap J*}_{(I\cup J)h'}\cN_{\ol{M}_{I\cap J}|\ol{M}})\sim &\\
& \sim \sum_{h'=g(h)} \alpha_{h'} \cdot_r \beta &
\eea
\end{proof}

\begin{lemma}
a) Given $\alpha\in Z_l(\ol{M}_J)$ and $I\subset J$, the following
relation holds: $$\bar{\phi}^I_{J*}(\alpha)\sim \sum_{[g]\in
G_I/G_J} g_*(\alpha),$$
where $g:\ol{M}_J\to\ol{M}_{g(J)}$ is the canonical isomorphism of Definition 3.5.

b)  The morphism
\bea \begin{CD} Z^*(\ol{M}(h)) @>{e(h)}>> B^*(\ol{M})
\\  \alpha @>>> \sum_{[g]\in G/G_h} \alpha_{g(h)} \end{CD} \eea
factors through rational equivalence.

c) The following sequence is exact:
$$\begin{CD} A^*(\ol{M}(h))@>{(j_*, -e(h))}>>
A^*(\ol{M})\oplus B^{' *}(\ol{M}) @>>> 
B^*(\ol{M})\to 0 \end{CD},$$
where $B^{' *}(\ol{M}):=\bigoplus_{l=0}^{\dim
(\ol{M})} ( \oplus_{I\not=\emptyset }
Z^{l-\codim_{\ol{M}}\ol{M}_{I}}(\ol{M}_{I}) / \sim )$. 
\end{lemma}
\begin{proof}
 
 Consider a cycle $\alpha=[V] \in Z_l(\ol{M}(h))$. Recall the corresponding
cycle $\alpha_h\in Z_l(\ol{M}_h)$ introduced for Definition 3.6.
   Assume that there is a sequence of distinct elements $h_1, ...h_m$ in the
$G$--orbit $Gh$, such that $V \subseteq \ol{M}(h_1...h_m)$ but $V
\not\subset \ol{M}(h_1...h_{m+1})$ for any $h_{m+1}\in Gh$ distinct from the
previous. The morphism $\ol{M}_{h}(h_2...h_m)\to \ol{M}(h h_2...h_m)$  is
$m:1$, which leads to the following decomposition: $$\alpha_h
=\frac{1}{m}(\alpha_h^1+...+\alpha_h^m)$$ (possibly $\alpha_h^i=\alpha_h^j$).
The same decomposition holds for $\alpha_{g(h)}$, where by convention,
$\alpha_{g(h)}^i:=g_*(\alpha_h^i)$ for all $[g]\in G/G_h$. 

 The following equality holds in $B^*(\ol{M})$:
\bean \sum_{[g]\in G/G_h} [\alpha_{g(h)}] = \sum_{[g]\in G/G_h}
[\alpha^i_{g(h)}] \eean for any $i\in \{ 1,...,m\}.$

 Indeed,  each  cycle $\alpha^i_{g_1(h)}\in
Z_l(\ol{M}_{g_1(h)}(g_2(h)...g_m(h))$ is  equivalent via $\sim$ to a sum of
cycles on  $\ol{M}_{ g_1(h)...g_m(h)}$, for all tuples of distinct $g_2,...,
g_m\in G/G_h$ allowable. Let $h_i:=g_i(h)$. Thus:
$$\sum_{h_1}\alpha^i_{h_1} \sim \sum_{h_1,...,h_m} \sum_{j=1}^{(m-1)!}
\alpha^{i, j}_{h_1...h_m} $$  where $\alpha^{i, j}_{h_1...h_m}$ are
the cycles of the $l$--dimensional components in $
\phi^{h_1^{-1}}_{h_1...h_2}(V^i_{h_1}),$ for the $(m-1)! :1$ map 
$\phi^{h_1}_{h_1...h_m}: \ol{M}_{h_1...h_m}\to \ol{M}_{h_1}.$
The multiplicity of each component $\alpha^{i, j}_{h_1...h_m}$ in
$\alpha^i_{h_1}$ is $1$ since there are $(m-1)!$ permutations of
$\{h_2,...,h_m\}$. On the other hand, the \'etale structure of the $G$--
network implies that for each allowable set $I=\{h_1,...,h_m\}\subset Gh$, the
preimages $\bigcup_{h_1}(\phi^{h_1}_{I})^{-1}(V^i_{h_1})$ and $\phi_I^{-1}(V)$
have the same $l$--dimensional components, with same multiplicities. Thus we
have shown that $\sum_{h_1}\alpha^i_{h_1}$ does not depend on $i$, which
proves the claim.

Claims a) and b) follow from the previous relation. Indeed, given
$\alpha\in Z_l(\ol{M}_J)$, it is enough to notice that $(\phi^I_{J *}\alpha
)_{h }$ contains $\alpha$ as a summand of its components for any $h\in
J\setminus I$, and then apply $\sim$ and relation (3.5) to obtain a).

 For b), consider  a cycle  $\alpha \in Z_l(\ol{M}(h))$ rationally
 equivalent to 0. We may assume for our purposes that $\alpha=[div(r)]$ for  a rational
 function $r$ on a $(l+1)$--dimensional scheme $W\subset \ol{M}(h)$. If we
 denote by $r_h$ the restriction of $r\circ \phi_h$ to the $(l+1)$--dimensional
components  of $\phi^{-1}_h(W)$, then $r_h^{-1}(0)$
  may capture only some of the $l$--dimensional components of
 $\phi^{-1}_h(r^{-1}(0))$. However, via relation (3.5),
$$ \sum_h [r^{-1}(0)]_h \sim \sum_h  [r_h^{-1}(0)]$$
and the same is true for $\infty$. Since $[div(r_h)]= 0$ in
$A_l(\ol{M}_h)$, it follows that $\sum_h [div(r)]_h = 0$ in
$B(\ol{M})$. This finishes the proof of b). c) is a direct consequence of b)
and the definition of $B^*$.  \end{proof}

The usual Chow ring of $\ol{M}$ may be recovered as invariant subring of
$B^*(\ol{M})$:

 \begin{lemma}
The action of $G$ on $\cP$ induces a canonical action on $B^*(\ol{M})$. Then
  $$A^*(\ol{M})= B^*(\ol{M})^G.$$ \end{lemma}
\begin{proof}
 
 By definition, the image of the natural morphism of
$A^*(\ol{M})\to B^*(\ol{M})$ is clearly included in
$B^*(\ol{M})^G$. A right inverse of this morphism is given
by the global  push--forward morphism $\Phi :B^*(\ol{M})\to
A^*(\ol{M})$. For $\alpha\in A_l(\ol{M}_{I})$, let $$\Phi (\alpha):=k_{|I|}
\bar{\phi}_{I*}\alpha,$$ where $k_I:=\frac{|G_I|}{|G|}$ is
chosen such that $\Phi$ is compatible with the equivalence relation in
$B^*(\ol{M})$.  $\Phi$ plays the role of Reynolds operator: it is
clearly invariant with respect to the $G$--action, and Lemma 3.13, a) implies
that $\Phi(\alpha)=\alpha$ for any $\alpha\in B^*(\ol{M})$. \end{proof}

%%incep demonstratia structura inel
%ma intorc la vechii indici
%atentie I in $P_{\geq d-k}$ ?

\subsection*{The extended Chow rings of the moduli spaces
$\ol{M}^k$ and their substrata} are computed here inductively after
the order of blow--ups. At step $k+1$, and for fixed $I\subset\cP$, the ring
$B^{*}(\ol{M}_I^{k+1}, k+1)$ is expressed as an algebra over
$B^{*}(\ol{M}_I^k, k)$.  We recall that $B^{*}(\ol{M}_I^k, k)$  is the
extended Chow ring of the network generated by all normalizations of the exceptional divisors of
$\ol{M}_I^k$ from previous blow--up steps. An intermediate algebra is
$B^{*}(\ol{M}_I^k, k+1)$, whose network has the $(k+1)$--th blow--up loci as
additional generators. Lemma 3.15 expresses $B^{*}(\ol{M}_I^k, k+1)$ as an
algebra over $B^{*}(\ol{M}_I^k, k)$ and Lemma 3.19 describes the algebra
$B^{*}(\ol{M}_I^{k+1}, k+1)$ over $B^{*}(\ol{M}_I^k, k+1)$. For this we
remember that modulo an equivalence relation $\sim$, $B^{*}(\ol{M}_I^{k+1},
k+1)$ is made by summands $B^{*+|J\setminus I|}(\ol{M}_J^{k+1}, k)$ for all
allowable sets $J$ such that $I\subseteq J\subset\cP$ and $J\setminus
I\subset \cP_{ d-k}$; similarly $B^{*}(\ol{M}_I^k, k+1)$  is made of
$B^{*+|J\setminus I|e}(\ol{M}_J^k, k)$ for the same set of $J$--s.  Lemma 3.16
identifies generators for the modules $B^{*}(\ol{M}_J^{k+1}, k)$ over
$B^{*}(\ol{M}_J^{k}, k)$, while a suitable pullback compatible with $\sim$ is
written in Lemma 3.17. The calculations are hardly surprising -- they
essentially transpose the usual Chow ring computations of \cite{fulton3},
\cite{keel} into our context. 

 A thread of induction goes through Lemmas 3.15 --3.19: we assume them
true at all previous blow--up steps and for all $J\supset I$. The Lemmas are
trivial at step 1, and for $J$ large enough that there is no element
$h\in\cP_{\geq d-k}\setminus J$ compatible with $J$, Lemma 3.4 applies since
 then $\ol{M}^{k+1}_J=\PP_J$.

\begin{notation}
We will employ the same notation $T_h:=-[\ol{M}^l_{Ih}]\in A^0(\ol{M}^l_{Ih})$,
for all $l> d-|h|$. The first Chern class of the normal line
bundle $\cN_{\ol{M}^l_{Ih}|\ol{M}^l_{I}}$ will be denoted by $-\tau_h$. 

We will omit the superscript $l$ for morphisms between $l$--th intermediate
spaces, as it will be visible in  the superscripts of the domains and
codomains. In particular, the blow--down morphism $\ol{M}^{k+1}_I\to \ol{M}^{k}_I$ will be denoted by $f_I$.  \end{notation}

\begin{lemma}
Assume that the pullback morphisms $\bar{\phi}^{I *}_J:  B^*(\ol{M}^k_I, k) \to
B^*(\ol{M}_J^k, k)$ are surjective for any $I \subset J$ such that
$J\setminus I \subset \cP_{\geq d-k}$. Then the  following $\QQ$--algebra isomorphism
holds: $$B^*(\ol{M}^k_I, k+1)  \cong B^*(\ol{M}^k_I, k)[\{y_h\}_{h\in
\cP_{d-k}\setminus I}] / \cI,$$ where the ideal $\cI$ is generated by:
\begin{enumerate}  \item $y_h^2=y_ha $, where
$a\in B^*(\ol{M}^k_I, k)$ is a class whose pullback to each $\ol{M}^k_{Ih}$ is
the first Chern class of the normal bundle $\cN_{\ol{M}^k_{h} | \ol{M}^k}$;
\item $y_hy_{h'}$ unless $h\cap h'=h, h' $ or $\emptyset$; 
\item $\prod_{h\in J\setminus I} y_h  \Ker \bar{\phi}^{I *}_J$, for any
$ J$ as above;  \item $(\sum_{h\in
\cP_{ d-k}\setminus I, h\cap h'=\emptyset } y_h) T_{h'}= -[\ol{M}_{Ih'}^k(h)]$
for any $h'\in \cP_{> d-k}\cup I$, where
$\ol{M}_{Ih'}^k(h)\hookrightarrow\ol{M}_{Ih'}^k$ is the image of
$\ol{M}_{Ihh'}^k$ in $\ol{M}_{Ih'}^k$, for any choice of $h\in
\cP_{ d-k}\setminus I$ such that $h\cap h'=\emptyset$. 
Here by  convention the formula also holds for $h'=\emptyset$, where 
$T_{\emptyset}$ is the unit element in $B^*(\ol{M}^k_I, k+1) $. \end{enumerate}

\end{lemma}

\begin{proof}
A morphism $F: B^*(\ol{M}^k_I, k+1)  \to  B^*(\ol{M}^k_I, k)[\{y_h\}_{h\in
\cP_{ d-k}\setminus I}]$ is defined on the  generators via the isomorphism
$$B^j(\ol{M}^k_I, k+1)\cong \bigoplus_{J\supseteq I, J\setminus I\subset
\cP_{ d-k}} B^{j-|J\setminus I|e}(\ol{M}_J^k, k) /\sim, $$ where
$e=\codim_{\ol{M}^k}\ol{M}^k_h=(n+1)(d-k)+1$, and the equivalence relation
$\sim$ previously defined  holds for $h\in \cP_{ d-k}\setminus I$.

Given $J$ as above and  $\alpha_J\in B^l(\ol{M}^k_J, k)$, there exists $\tilde{\alpha}_J$ in $B^{l+|J\setminus I|e}(\ol{M}^k_I, k)$ such that
$\alpha_J=\bar{\phi}^{I *}_{J}  \tilde{\alpha}_J$. Define $F(\alpha_J):=
\tilde{\alpha}_J\prod_{h\in J\setminus I} y_h$. Relation (3) insures that the
definition does not depend on the choice of the class $\tilde{\alpha}_J$.
Moreover,  compatibility with $\sim$  reduces to relation (4) via the
surjectivity of $\bar{\phi}^{I *}_{Jh}$ for any $h\in \cP_{ d-k}$. Finally,
relation (1) makes $F$ into a morphism of rings. The inverse morphism is
identity on $ B^*(\ol{M}^k_I, k)$ and takes $y_h$ to the class
$[\ol{M}^k_{Ih}]\in  B^0(\ol{M}^k_{Ih}, k)$.

\end{proof}

\begin{lemma}
 The generators of the $B^*(\ol{M}^k_I, k)$--module $B^*( \ol{M}^{k+1}_I, k)$
are linear combinations of products $\prod_{h\in J\setminus I}
  T_h^{l_h}$, for all allowable $J\supseteq I $ such that $J\setminus I\subset
\cP_{ d-k}$,  and $0\leq l_h \leq e$. \end{lemma}

% conventie asupra a cine e elementul unitate.

\begin{proof}

 The usual open--closed exact sequence implies the surjectivity of the map between graded  modules:
$$f^*B^*(\ol{M}^k_I, k)\oplus B^*(\ol{M}^{k+1}_I(h), k) \to B^*(
\ol{M}^{k+1}_I, k),$$ for $h\in \cP_{d-k}\setminus I$. Any generator $\beta$ of
$j_{Ih *}^IB^*(\ol{M}^{k+1}_I(h), k)$ is an element of
$j_{Jh *}^JA^*(\ol{M}^{k+1}_J(h))$ for some nested set $J\supseteq I$ such that
$J\setminus I\subset \cP_{>d-k}$, therefore being invariant under the action
of the group $G_{J\setminus I}\subset S_d$ which leaves all elements of
$J\setminus I$ invariant. On the other hand, $j_{Jh *}^JB^*( \ol{M}^{k+1}_J(h),
k)$ is embedded in  $B^{'*}( \ol{M}^{k+1}_J, k+1)$, which is generated by 
$B^{*+e|K\setminus J|}( \ol{M}^{k+1}_K, k)$, for all $K\supset J$ such that
$K\setminus J\subset \cP_{ d-k}$.  By induction, generators of 
$B^{*}(\ol{M}^{k+1}_K, k)$ are written over $B^{*}(\ol{M}^{k}_K, k)$   as
polynomials in $\{ T_h\}_{h\in\cP_{ d-k}}$. It follows that $\beta$ has a
polynomial expression of variables $\{ y_h\}$ and $\{ T_h\}$ over
$B^*(\ol{M}^{k}_J, k)$, for all  $h\in \cP_{ d-k}$, the exponents of $T_h$ not
exceeding $e$ by Lemma 3.4. Let $a\in B^e( \ol{M}^{k}, k)$ be a class whose
pullback to $\ol{M}^k_h$  is $c_{top}(\cN_{\ol{M}^k_h|\ol{M}^k})$.  Invariance
to $G_{J\setminus I}$, combined with relations $$y_h^2=y_h a \mbox{ and } T_hy_h=T_h a$$
imply that $\beta$ has in fact a polynomial expression of variables $\{ T_h\}$
and $\{ y_J(h)= \sum_{h'}y_{h'}\}$, over $B^*(\ol{M}^{k}_J, k)$, where  $h\in
\cP_{ d-k}$ and the sum is taken after all $h'\in \cP_{ d-k}$ having the same
set of incidence relations with the elements of $J$ as $h$. This
concludes the proof, as $ y_J(h) \in B^*(\ol{M}_J^k, k)$. We note that $\beta
$ also exhibits invariance to $G_{J\setminus I}$ as a polynomial in $\{ T_h\}$, but this
condition does not concern us here. It will be useful for Corollary 3.18.

%atentie la injectivitatea lui $j_*$

\end{proof}

\begin{notation} Set $e:=(n+1)(d-k)+1$.
 We will denote by $P_{\ol{M}^k_h|\ol{M}^k}(t) \in B^*(\ol{M}_I^k,k)[t]$
a fixed polynomial \bea P_{\ol{M}^k_h|\ol{M}^k}(t)= t^e+\sum_{i=0}^{e-1}
a_{h,i} t^i \eea whose pullback to
$B^*(\ol{M}^k_{Ih},k)[t]$ is $P_h$, the polynomial of Lemma 3.4, and whose free term is $y_h$.

Define $Q_h(t):=\frac{P_{\ol{M}_h^k |\ol{M}^k }(0)-P_{\ol{M}_h^k |\ol{M}^k }(t)}{t}.$

\end{notation}
$\tau_h$ will denote the first Chern class of the conormal line bundle
$\cN^{\vee}_{\ol{M}^{k+1}_h|\ol{M}^{k+1}}$, as well as its pull--backs to any
$\ol{M}^{k+1}_{Ih}$.  Recall that $\tau_h^i=-T_h^{i+1}$ by the
multiplication rule in $B^*(\ol{M}^{k+1}, k+1)$. Thus when $h\in I$, $\tau_h$
and $T_h$ will designate the same class on $\ol{M}^{k+1}_I$. 

 %atentie la semn

\begin{lemma}
  Assume that the pullback morphisms $\bar{\phi}^{I *}_J:  B^*(\ol{M}^k_I, k)
\to B^*(\ol{M}_J^k, k)$ are surjective for any $I \subset J$ such that
$J\setminus I \subset \cP_{ d-k}$.  
Then the morphisms of $\QQ$--vector spaces $ \bar{f}_J^*: B^*(\ol{M}^k_J, k) 
\to B^*(\ol{M}^{k+1}_{J}, k),$ $$ \bar{f}_J^*(\alpha):=
f_J^*(\alpha)\prod_{h\in J} Q_h( \tau_h)$$ for all allowable sets $J$
as above  add to a global pullback $$ \bar{f}^*_I:
B^*(\ol{M}^k_I, k+1)  \to B^*(\ol{M}^{k+1}_I, k+1).$$ \end{lemma}

\begin{proof}
  Due to the surjectivity of the maps $\bar{\phi}^{I *}_J$, compatibility
with the equivalence relation  $\sim$  reduces to the following fact: given
$J\supset I$ as above, and given $y:= [\ol{M}_{J}^k(h)] \in
A^0(\ol{M}_{J}^k(h))$ and $y_h:= [\ol{M}_{Jh}^k] \in A^0(\ol{M}_{Jh}^k) $ such
that $j^{J}_{Jh *} y = \sum_h y_h$ in  $ B^*(\ol{M}^k_I, k+1)  $,  then \bean  
f^{ *}_Jj^{J}_{Jh *} y  =  \sum_hQ_h( \tau_h) \eean  in  $B^*(\ol{M}^{k+1}_{I},
k+1)$.  This is proved by decreasing
induction on $|J|$. Indeed, $ f^{ *}_Jj^J_{Jh *} y  =  j^J_{Jh *} \alpha =
\sum_h \alpha_h $ for some class $\alpha \in A^0(\ol{M}_{J}^{k+1}(h))$.
Moreover by the previous lemma,  $j^J_{Jh *}\alpha $ may be written as a polynomial of
variables $T_h:=-[\ol{M}_{Jh}^{k+1}]\in A^*(\ol{M}_{Jh}^{k+1})$ over $
B^*(\ol{M}^k_J, k)  $, of degree no larger than $e$ in each variable.
Moreover, for dimension reasons, $T_h^e$ can only appear in $j^J_{Jh *}\alpha$ with a
numerical coefficient  $b_h$. In fact, $b_h=-1$, as can be seen via the
push--forward $f_{h*}: B^*(\ol{M}_{Jh}^{k+1}, k+1) \to
B^*(\ol{M}^{k+1}_{J}(^k_h), k+1) $ (as in the end of section 2), because $f_{h*}$  maps each
monomial $T_h^{l_h}$ into zero whenever $l_h<e$, and maps $T_h^e$ into
$y_h$. It remains to find the other coefficients.

 Pullback to $\ol{M}_{Jh}^{k+1}$ and the intersection formula of Lemma 3.9
yield, on the one hand:  $$ \bar{\phi}^{J *}_{Jh}j^J_{Jh *}\alpha=  
\bar{\phi}^{J *}_{Jh} f^{*}_Jj^J_{Jh *} y = f^{*}_{Jh}\bar{\phi}^{J *}_{Jh}j^J_{Jh
*} y = f^{*}_{Jh}(c_{e}(\cN_{\ol{M}_h^k |\ol{M}^k }) +   y_h(h'))$$ where
$y_{h}(h')= \bar{\phi}^{Jh}_{Jhh' *} \bar{\phi}^{Jh' *}_{Jhh'} y_{h'}$ is the class
of the blow--up locus in  $\ol{M}_{Jh}^{k}$.    As a
polynomial over $B^*(\ol{M}_{Jh}^{k+1}, k+1)$, $\bar{\phi}^{J *}_{Jh}j^J_{Jh
*}\alpha$ preserves the form of $j^J_{Jh *}\alpha$. Pullback from $\PP_h$ of relation in Lemma 3.4 yields   $$
c_{e}(\cN_{\ol{M}_h^k |\ol{M}^k })  = \tau_h Q(\tau_h ).$$  By
induction,   $f^{ *}_{Jh} y_h(h')=\sum_{h'} T_{h'}Q(T_{h'}),$ the sum being
after all $h'\in\cP_{ d-k}\setminus Jh$ compatible with $h$. Therefore in
$B^*(\ol{M}^{k+1}_J, k+1)$, $$j^J_{Jh *}\alpha = T_h Q(T_h )+\sum_{h'} T_{h'}Q(T_{h'}),$$
modulo $\Ker \bar{\phi}^*_h: B^*(\ol{M}^k, k)\to  B^*(\ol{M}^k_h, k)$. Applied
for all $h\in \cP_{ d-k}$, this  completely determines the coefficients
of $\alpha$ as a polynomial in $\{T_h\}_h$. This ends the compatibility check.

It remains to prove that $\bar{f}^*_I$ preserves the multiplication structure.
Given the definition  of $\bar{f}^*_I$, the only relevant case is that of
$\bar{f}^*_I(y_h^2)$, when $y_h:= [\ol{M}_{Ih}^k]\in A^0(\ol{M}_{Ih}^k)$. Then
$y_h^2= a y_h$, where $a\in  B^*(\ol{M}^k, k) $ is a class whose pullback to 
$B^*(\ol{M}^k_{h}, k)$ is $c_{e}(\cN_{\ol{M}_h^k |\ol{M}^k })$. Thus $
\bar{f}^*_I(y_h^2)= -T_h Q(T_h) a$. On the other hand, $ \bar{f}^*_I(y_h)
\bar{f}^*_I(y_h)=-T_h^2 Q(T_h)^2$, which equals the former expression since
$\tau_h Q(\tau_h)= c_{e}(\cN_{\ol{M}_h^k |\ol{M}^k })$ on $\ol{M}^k_{Ih}$.

\end{proof}

\begin{corollary}
 An element $\beta\in B^*(\ol{M}^{k+1}_I, k)$ is zero if and only if $T_h\beta
=0$ for all  $h\in \cP_{d-k}$ and $f_{I *}(\beta)=0$ in 
$B^*(\ol{M}^k_I, k)$. \end{corollary}

\begin{proof}

Recall that the map $f_I$ factors through $g_I: \ol{M}^{k+1}_I\to \PP_I$ and
$p_I: \PP_I \to \ol{M}^k_I$. On the weighted projective fibration $p_I: \PP_I
\to \ol{M}^k_I$ the result is trivial. Thus  we may  assume $g_{I
*}(\beta)=0$.

By  Lemma 3.16, an element $\beta\in B^*(\ol{M}^{k+1}_I, k)$ can be written
as a polynomial of variables $\{ T_h\}$, $h\in \cP_{ d-k}$ over $B^*(
\ol{M}^{k}_I, k)$ , the exponents of $T_h$ not exceeding $e$. This
can be refined to exponents strictly less than $e$, due to the
invariance of $\beta $ under some group $G_{J\setminus I}$, discussed in the proof of Lemma
3.16,  and to relation (3.6) in  Lemma 3.17.  
Consequently, $g_{I*}(\beta)$ is the term free of variables $\{ T_h\}_{h\not\in
I}$.

 Fix $h\not\in I$. The pullback $\bar{\phi}_{Ih}^{I*}\beta$ has the same
polynomial expression as $\beta$, but over $B^*(\ol{M}^k_{Ih},k)$. As 
$\ol{M}^{k+1}_{Ih}$ is a weighted projective fibration in the $h$--direction, the
annihilation of $T_h\beta$ implies that $\beta$, as a polynomial in $T_h$, has
all coefficients in  $\Ker \bar{\phi}_{Ih}^{I*}\subset B^*(\ol{M}^k_I, k)$ (see
Lemma 3.4). Inductively we reason that the coefficient of $\prod_{h\in
J\setminus I}T_h^{l_h}$ in $\beta$ is an element of  $\Ker \bar{\phi}_J^{I *}$, for
each allowable $J$ with $J\setminus I\subset \cP_{ d-k}$.  This concludes the
proof, as the relation $$\prod_{h\in J\setminus I}T_h^{l_h}\Ker
\bar{\phi}_J^{I *} =0$$ is inherent in the multiplicative structure of
$B^*(\ol{M}^{k+1}_I, k+1)$.

\end{proof}

\begin{lemma}

Under the same conditions as in the preceding lemmas, the  following
$\QQ$--algebra isomorphism holds: $$B^*(\ol{M}^{k+1}_I, k+1)  \cong
B^*(\ol{M}^k_I, k+1)[\{T_h\}_{h\in \cP_{ d-k}}] / \cJ,$$ where the ideal $\cJ$ is
generated by: \begin{enumerate}
\item $T_hT_{h'}$ unless $h\cap h'=h, h' $ or $\emptyset$;
\item $\prod_{h\in J\setminus I} T_h  \Ker \bar{\phi}^{I *}_J$, for any
$J\supseteq I$ such that  $J\setminus I\subset \cP_{ d-k}$; \item $P_{\ol{M}_h^k
|\ol{M}^k }(T_h) $, for all  $h\in \cP_{ d-k}$. \end{enumerate}
\end{lemma}

\begin{proof}

We denote  by $QB^*(\ol{M}^{k+1}_I, k+1)$ the  $ B^*(\ol{M}^k_I, k+1)$--
algebra  $B^*(\ol{M}^k_I, k+1)[\{T_h\}_{h\in \cP_{ d-k}}] / \cJ$. There is an
obvious morphism of algebras: $$\Phi : QB^*(\ol{M}^{k+1}_I, k+1) \to
B^*(\ol{M}^{k+1}_I, k+1),$$ sending each generator $T_h$ into the class 
$-[\ol{M}^{k+1}_{Ih}] \in  B^0(\ol{M}^{k+1}_{Ih}, k).$ Indeed, by
Lemma 3.17 and the multiplicative structure of $B^*(\ol{M}^{k+1}_I, k+1)$, all
generators of $\cJ$ map to zero in $B^*(\ol{M}^{k+1}_I, k+1)$. It remains to
construct the inverse of $\Phi$.

Via the surjective map $ f_I^*B^*(\ol{M}^k_I,k)\oplus
B^*(\ol{M}^{k+1}_I(h),k)\to B^*(\ol{M}^{k+1}_I,k)$,  the exact sequence of
Lemma 3.13, c) can be further refined to $$  BK^*(\ol{M}^{k+1}_I,k)
\to f_I^*B^*(\ol{M}_I^k,k)\oplus B^{'*}(\ol{M}^{k+1}_I,k+1)\to
B^*(\ol{M}^{k+1}_I,k+1) $$ where $BK^l(\ol{M}^{k+1}_I,k)=\{
\alpha \in B^l(\ol{M}^{k+1}_I(h),k)  |  j_{I*}\alpha \in f_I^*j_{I*}
B^*(\ol{M}^{k}_I(h),k)  \}$.  Let the rightmost morphism above be $\Psi$. There is a
natural morphism of $B^*(\ol{M}^k_I,k)$--modules  $$F:B^*(\ol{M}^k_I,k)\oplus
B^{' *}(\ol{M}^{k+1}_I,k+1)\to  QB^*(\ol{M}^{k+1}_I,k+1).$$ Indeed, $B^{'*}(\ol{M}^{k+1}_I,k+1)$ is
generated by the $B^*(\ol{M}^k_J,k)$--modules $B^*(\ol{M}^{k+1}_J,
k)$, for $J\setminus I\subset \cP_{ d-k}$, which are generated by monomials
$\prod_{h\in J}T_h^{l_h} $, and we may simply define $$F( \prod_{h\in
J}T_h^{l_h}):= \prod_{h\in J} T_h^{l_h}.$$ The image of $B^*(\ol{M}^k_I,k)$ in
$QB^*(\ol{M}^{k+1}_I,k+1)$ is naturally defined via the inclusion 
$B^*(\ol{M}^k_I,k) \hookrightarrow B^*(\ol{M}^k_I,k+1)$. The construction of
$F$ implies $\Phi\circ F=\Psi$. It remains to show that the image of
$BK^*(\ol{M}^{k+1}_I,k)$ through $F$ is zero.

%de mentionat inductia, primul pas de inductie;

 We recall that $j_{I*} B^*(\ol{M}^{k}_I(h),k)\hookrightarrow
B^{'*}(\ol{M}^{k}_I,k+1)$ is generated over
$B^*(\ol{M}^{k}_I,k)$ by classes $[y_h]$ with $h\in \cP_{ d-k}$. Thus the
elements in the image of the morphism
$$    BK^*(\ol{M}^{k+1}_I,k)
\to f_I^*B^*(\ol{M}_I^k,k)\oplus B^{'*}(\ol{M}^{k+1}_I,k+1)  $$
are linear combinations of  $(y_h, -\{Q(T_h)T_h\}_{h\in \cP_{ d-k}})$
for  $h\in \cP_{ d-k}$, and of elements in $(0, \Ker j_{I *})$ over 
$B^*(\ol{M}^k_I, k+1)$. The first are mapped to zero in $QB^*(\ol{M}^{k+1},
k+1)$ by relation (3) in the ring above, and the second type of elements are
also mapped into zero by the same proof as of Corollary 3.18, because of relation
(2).

\end{proof}

%strict images

 The following lemma is an adaptation to our situation of Lemma 5.2
 from \cite{fulton2}. It presents an inductive way of constructing the
 classes of the blow--up loci, by showing how such
 classes are transformed under the preceding steps of weighted blow--up.

  Fix $I, J\subset \cP$ such that $I\cap J=\emptyset$. Recall that
$\ol{M}^l_I(J)$ is the image of $\bar{\phi}^I_{I\cup J}:\ol{M}^l_{I\cup J}\to
\ol{M}^l_I$. Here we consider the case $J=\{ h'\}$, and $|h'|=d-k'<d-k$. 

\begin{lemma}
 The class of the strict transform $\ol{M}_{I}^{k+1}(h')$ of $\ol{M}_{I}^k(h')$
in $B^*(\ol{M}^{k+1}_I, k+1)$ is given by $$[\ol{M}_{I}^{k+1}(h')] =
f_I^*([\ol{M}_{I}^k(h')])+ \sum_{h\in \cP_{ d-k}}
(P^{h}_{\ol{M}_{h'}^k|\ol{M}^k}(T_h) - P^{h}_{\ol{M}_{h'}^k|\ol{M}^k}(0))$$
 where $h'$ is a subset of $h$ and
$P^{h}_{\ol{M}_{h'}^k|\ol{M}^k}(t)$ is a polynomial over $B^*(\ol{M}^k_I,k)$
whose pull--back to $\ol{M}_{Ih}^k$ is
$P_{\ol{M}_{Ih}^k|\ol{M}^k_I}(t)/P_{\ol{M}_{Ih}^k|\ol{M}_{Ih'}^k}(t)$.

Note that this polynomial does not depend on the choice of
$h'\in\cP_{d-k'}$, as long as $h'\subset h$, and neither do the classes
$[\ol{M}_{I}^k(h')]$ and $[\ol{M}_{I}^{k+1}(h')]$, as the group $G_h$
 permutes the elements of $h'\in\cP_{d-k'}$.

\end{lemma}

\begin{proof}
By Corollary 3.18, it is enough to show that the above formula holds
 when we apply $f_{I *}$ and $\bar{\phi}_{Ih}^{I*}$ to both
 sides. The first part is a consequence of the equality
 $$f_{I*}(P^{h}_{\ol{M}_{h'}^k|\ol{M}^k}(T_h) -
P^{h}_{\ol{M}_{h'}^k|\ol{M}^k}(0))=0,$$ since the polynomial above has no free
terms in $T_h$. The second part comes out of descending induction on $|I|$.

 The class of the preimage $(\bar{\phi}_{Ih}^{I})^{-1}(\ol{M}_{I}^k(h'))$ is
the sum of two distinct components: a constant multiple of $[\ol{M}_{Ih}^k]$,
corresponding to the set $\{h'\in \cP_{d-k'}  | h'\subset h \}$, and
$[\ol{M}_{Ih}^k(h')]$ corresponding to the set $\{h'\in \cP_{d-k'}  |
{h'}\cap h=\emptyset \}$. 
 Since in the last case, the preimages of $\ol{M}_{Ih'}^k$ and 
$\ol{M}_{Ih}^k$ intersect transversely in an \'etale cover of $\ol{M}^k_I$, a
standard dimension argument shows that
$f_{Ih}^*[\ol{M}_{Ih}^{k}(h')]=[\ol{M}_{Ih}^{k+1}(h')]$.  For the former case, 
the induction hypothesis on $\ol{M}_{Ih}^{k+1}$ states that
$$[\ol{M}_{Ih}^{k+1}(h')]=g_{Ih}^*[\PP_h(h')]+
\sum_{j\in \cP_{ d-k}} (P^{j}_{\ol{M}_{j'}^k|\ol{M}^k}(T_j) -
P^{j}_{\ol{M}_{j'}^k|\ol{M}^k}(0)) $$ where the sum is taken after all $j$
with $\{j, h\}$ allowable, and  $j'\subset j$ is fixed for each $j$.
Here $\PP_h(h')$ is the weighted projective fibration of
$\cN_{\ol{M}_{h}^k|\ol{M}^k_{h'}}$, pulled back to $\ol{M}^k_{Ih}$. It is
locally the quotient $\bigcup_{h'\subset h}
\PP(\cN_{\ol{M}_{Ih}^k(\bar{t})|\ol{M}^k_{Ih'}(\bar{t})})/G_{Ih}$. The
morphism $f_{Ih}: \ol{M}_{Ih}^{k+1}\to \ol{M}_{Ih}^{k}$ factors
through  $g_{Ih}: \ol{M}_{Ih}^{k+1}\to
\PP_h\times_{\ol{M}^k_h}\ol{M}^k_{Ih}$ and $\pi_h:
\PP_h\times_{\ol{M}^k_h}\ol{M}^k_{Ih}\to  \ol{M}_{Ih}^{k}$, while
$g_{Ih}(\ol{M}_{Ih}^{k+1}(h'))=\PP_h(h')$. 

 It remains to show that the class of $\PP_h(h')$ in
$\PP_h\times_{\ol{M}^k_h}\ol{M}^k_{Ih}$ is

$$ [\PP_h(h')]=P^{h}_{\ol{M}_{h'}^k|\ol{M}^k}(T_h)
- P^{h}_{\ol{M}_{h'}^k|\ol{M}^k}(0). $$

   This is proven by the same methods as in Lemma 3.2: Let $	\tau_h$ be the
class of the line bundle $\cO_{\PP_h}(1)$. After writing the class
$[\PP_h(h')]$ as    
$$[\PP_h(h')]=\sum_{i=0}^{e-1} b_i\tau_h^i,$$ where
$e=\codim_{\ol{M}^k}\ol{M}^k_h$,   and considering
$\pi_{h*}([\PP_h(h')]\tau_h^l)$ for all $l$, we obtain that
the sum $\sum_i b_i t^i$ satisfies:  $$(\sum_i b_i t^i)(\sum_i
\pi_{h*}(\tau_h^i)t^{-i-1})=\sum_i
\pi_{h*}([\PP_h(h')]\tau_h^i)t^{-i-1}.$$   The right side is
the equivariant Segre polynomial of $\cN_{\ol{M}^k_{h}|\ol{M}^k_{h'}}$ in the
sense of Lemma 3.2 , i.e. the inverse of $P_{\ol{M}^k_{h}|\ol{M}^k_{h'}}$, and
$\pi_{h*}(\tau_h^i)t^{-i-1}$ is the inverse of $P_{\ol{M}^k_{h}|\ol{M}^k}$. This finishes the proof of the
lemma. \end{proof} 

It remains to compute the polynomials $P_{\ol{M}^k_{Kh}|\ol{M}^k_K}$, for a
nested set $K$ and $h\in\cP_{ d-k}\setminus K$. The single important case is
when $K$ is a partial partition. We may split $K$ into two partial partitions
$I$ and $J$, one containing subsets of $h$, the other containing sets disjoint
from $h$.

\begin{notation}
 Consider two partial partitions $I$ and $J$ of $M_0$, such that
 $IJ:= I\bigcup J$ is still a partial partition.
Let $h_0\in \cP_{d-k_0}$ such that $h_0\supseteq \cup_{h\in I} h$. To these and
 to any number $k$ such that $k_0 \leq k\leq d-|\cup_{h\in I} h|$, we associate
  the polynomial $P^k_{I,J}(\{t_h\}_{h \supseteq h_0})$ of variables
 $\{t_h\}_{h\supseteq h_0}$ over $\QQ[H,\psi]$, defined as follows:
%\bea & P^k_{I,J} (\{t_h\}_{h\supseteq h_0})=(\psi+\sum_{h\supseteq h_0}
% t_h)^{|I|-1}\cdot & \\
%&  \cdot\prod_{j=k-l_J+1}^{d-l_{I\cup J}} (H+j\psi +\sum_{h\supseteq
%h_0} (|h|+l_{J_{h}}-d+j)T_h)^{n+1}- &\\ &- (H+(k-l_J)\psi
%+\sum_{h\supseteq h_0}(|h|-l_{J\backslash J_h}-d+k)
% T_h)^{(d-l_I-k)(n+1)} &\eea

\bean P^k_{I,J}(\{t_h\}_{h \supseteq
h_0})=    \eean \bea
= \psi_{h_0}^{|I|-1} (\prod_{j=1+|h_0|-(d-k)}^{|h_0\setminus
\cup_{h\in I}h|}(H_{h_0}+j\psi_{h_0})^{n+1}-
H_{h_0}^{(n+1)(d-k-|\cup_{h\in I}h|)}), \eea where \bea
\psi_{h_0}=\psi+\sum_{h\supseteq h_0} t_h, \eea \bea
H_{h_0}=H+(d-|h_0\cup(\cup_{h\in J}h)|)\psi+\sum_{h\supset
h_0}|h\backslash (h_0\cup(\cup_{h'\in J}h'))|t_h.\eea We denote
by \bean P^k_{I,J}(0):=\psi^{|I|-1} \prod_{j=1}^{d-k-|\cup_{h\in
I}h|}(H+(k-|\cup_{h\in J}h|+j)\psi)^{n+1}. \eean
\end{notation}

 $P^k_{I,J}$ admits the
following geometrical
 interpretation:
%the normal bundle
\begin{lemma}
 Let $h\in \cP_{\geq d-k}$ be
such that $h\supset \bigcup_{h'\in I}h'$  and ${h}\bigcap h'=\emptyset $  for all $h'\in J$.  Then the 
following relation holds:
$$P_{\ol{M}^{k}_{Jh} |\ol{M}^{k}_{IJ}}(t)=P^k_{I,J}(\{ T_{h'}\}_{h'\supset h}, t). $$
The class of $[ \ol{M}_{IJ}^k(h)]$ in $[ \ol{M}_{IJ}^k]$ is
$$[ \ol{M}_{IJ}^k(h)]=\left.\sum_{h_0} P^k_{I,J}(\{T_{h'}\}_{h' \supset h_0},
t_{h_0})\right|^{t_{h_0}=T_{h_0}}_{t_{h_0}=0}+ P^k_{I,J}(0),$$
where the sum is taken after all $h_0\in \cP_{d-k_0}$ with $k_0<k$ and
such that $h_0\supset \cup_{h'\in I} h'$ and $I\cup J\cup\{h_0\}$
is allowable.

\end{lemma}
\begin{proof}
%Consider the commutative diagram
%$$ \diagram
%\ol{M}^{k+1}_{h;I,J} \rto^{} \dto_{f_{h}(I,J)} & \ol{M}^{k+1}_{IJ} \dto^{} \\
%\ol{M}^{k}_{Jh}  \rto_{} & \ol{M}^{k}_{IJ}
%\enddiagram. $$
  The map $f_{IJh}:\ol{M}^{k+1}_{IJh} \to \ol{M}^{k}_{Jh}$ factors through
$\PP_{h;IJ}\to\ol{M}^{k}_{Jh}$, where $\PP_{h;IJ}$ is the weighted projective
fibration of the normal bundle $\cN_{\ol{M}^{k}_{Jh}|\ol{M}^{k}_{IJ}}$. 
%tthis notion to be defined somewhere

It was remarked in section 2 that $\PP_{h;IJ}\cong \ol{M}_{I}^{1,
|h|}\times_{\PP^n}\ol{M}^{k}_{Jh}$,  where $\ol{M}_{I}^{1,
|h|}$ is  the space of 1--stable, $I$-- type, degree $|h|$
rational maps, and the normalization of a closed stratum in $\ol{M}^{1,
|h|}:=\ol{M}_{0,1}(\PP^n, |h|, 1)$. Its associated polynomial
\bean  P_{\ol{M}_{I}^{1,
|h|}|\ol{M}^{1,|h|}}=   
\psi_h^{|I|-1}\prod_{i=1}^{|h|-l_I}(H_h+i\psi_h)^{n+1}, \eean where $H_h$ is the
pull--back of the hyperplane section in $\PP^n$ and $\psi_h$ is the cotangent
class.

On the other hand, there is a natural birational morphism
$$\ol{M}_{0, 1+|Jh|}(\PP^n, d-|h\cup (\cup_{h'\in J}h')|)\to \ol{M}^{k}_{Jh}$$
making $ \ol{M}^{k}_{Jh}$ into an intermediate moduli space
of $\ol{M}_{0, 1+|Jh|}(\PP^n, d-|h\cup (\cup_{h'\in J}h')|)$.
 The pullback $H_h$ of the hyperplane class in $\PP^n$ by the evaluation
morphism  $e_{h}$: $\ol{M}^{k}_{Jh}\to \PP^n$ at the point marked by $h$ is
given by the following expression \bea H_h:= H+ (d-|h\cup (\cup_{h'\in
J}h')|)\psi +\sum_{h\subset {h''}} |h''\setminus (h\cup (\cup_{h'\in
J}h'))| T_{h''},  \eea conform Theorem 1 in
\cite{lee-pandharipande}).
  This class is clearly pulled--back from
$\ol{M}^{k}_{Jh}$ and its pullback to $\PP_{h;IJ}$ is the class denoted
also by $H_h$.

 The map 
$\ol{M}^{k+1}_{hJ} \to \ol{M}^{k+1}_{IJ}$ pulls back to the
gluing map $f$: \bea \begin{CD} \ol{M}_{0, 1+|Jh|}(\PP^n,
d-|h\cup (\cup_{h'\in J}h')|)\times_{\PP^n}\ol{M}_{0,1+|I|}(\PP^n, |h\setminus (\cup_{h'\in I}h')|) \\ @V{f}VV \\
\ol{M}_{0,1+|IJ|}(\PP^n, d-|\cup_{h'\in I\cup J}h'|)\end{CD}.\eea
  Let $\pi_1, \pi_2$ be the projections  from the domain of $f$
 onto $ \ol{M}_{0, 1+|Jh|}(\PP^n,
d-|h\cup (\cup_{h'\in J}h')|)$ and $\ol{M}_{0,1+|I|}(\PP^n,
|h\setminus (\cup_{h'\in I}h')| )$. Let $\psi_i$ be the cotangent
line class coming from the $i$--th marked point. The following
relation is also part of Theorem 1 in \cite{lee-pandharipande}:
$$ -\pi_1^*\psi_{h} -\pi_1^*\psi_1 = \sum_{h\subset h'} f^*T_{h'}.$$
The following are also known to be true: $$ \pi_2^*\psi_1 +
\pi_1^*\psi_{h} = f^* T_{h} \mbox{ and } \pi_1^*\psi_1 =
f^*\psi_1.$$ Putting these together:
 \bea \pi_2^*\psi_1 = f^*\psi_1 + \sum_{h' \supseteq {h}} f^*T_{h'}.\eea
  This is clearly pullback of the class $\psi_h$ from
  $\PP_{h;IJ}$.
 Thus $H_h$ and $\psi_h$ have been written as linear combinations of $H$,
$\psi$ and boundary divisors. This, together with formula (3.9), determine the
polynomial $P_{\ol{M}^k_{hJ}|\ol{M}^{k}_{IJ}}(T_h)$.
The  term $H_h^{n+1}$ was introduced in expression (3.7) for
division purposes in the case when $I=\emptyset$. It is clearly zero since it
comes form a null class in $\PP^n$.

 The second part follows after applying Lemma 3.20 through successive blow--up
steps.

\end{proof}
 The next proposition gives a preliminary formula for
$B^*(\ol{M}_{0,1}(\PP^n,d)):=B(\ol{M}^d,d)$:

\begin{proposition}
 The $\QQ$--algebra $B^*(\ol{M}_{0,1}(\PP^n,d))$ is
 generated over $\QQ$ by the divisor classes $H$, $\psi$ and $\{T_h\}_h$ for
all strict   subsets $h$ of $M=\{1,...,d\}$. Its ideal of relations
is generated by: \begin{enumerate}
 \item $ H^{n+1};$
\item $T_hT_{h'}$ for all $h,h'$ such that $\{h, h'\}$ is not allowable,
i.e. if $h\bigcap {h'}\not= \emptyset, h$ or ${h'}$.
\item $ \prod_{h\in I} T_h \psi^{|I|-1}\prod_{i=1}^{d-l_I}(H+i\psi)^{n+1}$
for any partial partition $I$ of $M$;
\item $ \prod_{h\in I\cup J} T_{h}\left(\sum_{h_0\in \cP_{\geq d-k}(I,J)}\left.
P^k_{I,J}(\{T_{h}\}_{h\supset h_0},
t_{h_0})\right|^{t_{h_0}=T_{h_0}}_{t_{h_0}=0}+ P^k_{I,J}(0)\right) $,
 for any
partial partitions $I$ and $J$ of $M$ such that  $I\bigcup J$ is still a
partial partition. Here $\cP_{\geq d-k}(I,J)$ is the set of all subsets $h_0$ of $M$
with ${|h_0|\geq d-k}$ and such that $h_0\supseteq \bigcup_{h \in I} h$
and $h_0\bigcap {h}=\emptyset$  whenever $|h_0|= d-k$ and $h\in
J$.

\end{enumerate}
\end{proposition}

\begin{proof}
The ring $B^*(\ol{M}_{0,1}(\PP^n,d))$ is constructed inductively.
 Let $K$ be any partial partition of $M$, and let $k\leq min_{h\in K}\{ d-|h|\}$.
By successive applycations of  Lemmas 3.15 and 3.19 we know that
$$B^*(\ol{M}^k_K, k)=\QQ[H, \psi,  \{T_h\}_{|h|\geq d-k}]/\cI_K$$ for some ideal
$\cI_K$. In particular, we may write $$B^*(\ol{M}_{0,1}(\PP^n,d))=\QQ[H,\psi, 
\{T_h\}_h]/\cI$$ where $h$ ranges over all strict subsets of $M$. From relation
(2) in Lemma 3.19 we conclude that for any element  $\alpha\in \cI$,  there
exists a maximal set $K$ and number $k$ such that  $\alpha \in  \prod_{h\in K}
T_h \cI_K.$    When $K=\emptyset$ and $k=1$ we obtain elements $H^{n+1}$ and 
$\psi^{-1}\prod_{i=1}^d(H+i\psi)^{n+1}$.  Elements (3) in $\cI$ arise by
Lemma 3.4  when $K\not=\emptyset$ and  $k=1$. Elements (2) in $\cI$ are
natural compatibility conditions among divisors.  For $K$ and $k>1$ fixed as
above, the blow--up locus of  $\ol{M}^k_K$ at step $k$ has a number of distinct
components, indexed by the choices of subsets $I\subset K$ such that
$\sum_{h\in I}|h|\leq d-k$. For such an $I$, and $J$ such  that $J\setminus (K\backslash I)\subset \cP_{>d-k}$, consider
$h\in\cP_{ d-k}$ such that $h\supseteq \bigcup_{h'\in I} h'$ and $h\bigcap h'=\emptyset$ for all $h'\in J$.  Putting together
 relations (4) of Lemma 3.15 with relation (3) in Lemma 3.19, plus the
polynomial calculations of Lemma 3.21, we obtain elements (4) of $\cI$. 

\end{proof}

As is the case with successive blow--ups, many of the relations coming from
intermediate steps are superfluous. The next theorem presents a simplified
version of the $B^*(\ol{M}_{0,1}(\PP^n,d))$ ring structure. The simplification
calculations form the Appendix. 

\begin{theorem}
 Consider the $\QQ$--algebra $B^*(\ol{M}_{0,1}(\PP^n,d);\QQ)$
 generated by the divisor classes $H$, $\psi$ and $\{T_h\}_h$ for all proper
 subsets $h$ of $M=\{1,...,d\}$, modulo the ideal generated by:
\begin{enumerate}
 \item $ H^{n+1};$
\item $T_hT_{h'}$ for all $h,h'$ such that $h\bigcap {h'}\not= \emptyset, h$ or ${h'}$;
\item $T_hT_{h'}(\psi+\sum_{h''\supseteq h\cup h'} T_{h''})$ for all $h$, $h'$ nonempty and disjoint;
\item $T_h(\sum_{h'\not\subseteq h} P^{d-1}_{\emptyset,
h}(\{T_{h''}\}_{h''\supset h'}, t_{h'}) \left.
\right|^{t_{h'}=T_{h'}}_{t_{h'}=0}+P^{d-1}_{\emptyset, h}(0))$ for all
$h\subset M$, where $T_{\emptyset}:=1$ and \bea &P^{d-1}_{\emptyset,
h}(\{t_{h''}\}_{h''\supseteq h'}):=(\psi+\sum_{h''\supset h'} t_{h''})^{-1}
[(H+d\psi+\sum_{h''\supseteq h'}|h''| t_{h''})^{n+1}-&\\
&-(H+(d-|h\cup h'|)\psi+\sum_{h''\supset h'}(|h''\setminus (h\cup h')|) t_{h''})^{n+1})]&\eea

\end{enumerate}

 The group of symmetries $S_d$ has a natural action on the set of proper subsets of $M$. 

$A^*(\ol{M}_{0,1}(\PP^n,d);\QQ)$ is the ring of
invariants of the ring $B^*(\ol{M}_{0,1}(\PP^n,d);\QQ)$ under the
induced
action.
\end{theorem}

 %d=2
\begin{example}
Here we describe in detail $ \ol{M}_{0,1} (\PP^n , 2) $ and its Chow ring.
In this case $M=\{1,2\}$. There is only one intermediate space $\ol{M}_{0,1} (\PP^n , 2, 1)$, with a 
 blow--up locus $\ol{M}^1(1)\subset \ol{M}_{0,1} (\PP^n , 2, 1)$ and  two
copies of  its normalization,  $\ol{M}^1_{\{1\} }$ and $\ol{M}^1_{\{2\} }$.

  In accord with our considerations in Lemma 3.3, $\ol{M}^1_{\{i\} }
\cong \PP^n \times \PP^n $ for $i\in \{1,2\}$.  Outside its diagonal
$\Delta \subset \PP^n\times\PP^n$, we may think of $\ol{M}^1_{\{i\} }$ as
parametrizing lines in $\PP^n$ with two marked points $q_1$ and $q_2$
on them: $q_1$ stands for the marked point of $\ol{M}_{0,1}(\PP^n,2)$,
$q_2$ is identified by the set $\{i\}$. The diagonal
is isomorphic to  $\ol{M}^1_{\{1\},\{2\} }$, parametrizing rational curves
with $3$ marked points, contracted to a point in $\PP^n$.  Thus the
two projections $\pi_1,\pi_2:\PP^n\times\PP^n\to \PP^n$ stand for
evaluation maps, the two corresponding pull--backs of the hyperplane
class in $\PP^n$ are $H$ and $H'=H+\psi$, and the class of
$\ol{M}^1_{\{1\},\{2\}}\cong \Delta$ in $\ol{M}^1_{\{i\} }
\cong \PP^n \times \PP^n $ is $$\sum_{k=0}^{n} H^k H'^{n-k}
=\frac{(H+\psi)^{n+1}}{\psi}.$$ On the other hand, by Lemmas 3.3 and
3.21, the image of  $\ol{M}^1_{\{i\} }$ in $\ol{M}^1$ has class
$$[\ol{M}^1_{\{i\} }]=\frac{(H+2\psi)^{n+1}}{\psi},$$
$P_{\ol{M}^1_{\{i\} }|\ol{M}^1}(t)$ is the polynomial $$P(t):=\frac{(H+2\psi
+t)^{n+1}-(H+\psi)^{n+1}}{\psi +t}.$$  The intersection formula (3.3) on
$\ol{M}^1_{\{i\} }$ reads:   $$\frac{(H+2\psi)^{n+1}}{\psi} =
\left.\frac{(H+\psi)^{n+1}}{\psi} + \frac{(H+2\psi
+t)^{n+1}-(H+\psi)^{n+1}}{\psi +t} \right|_{t=0}.$$   The exceptional divisor in $\ol{M}^2$ is 
$\ol{M}^2_{\{i\} }\cong Bl_{\Delta}(\PP^n\times\PP^n)
\times_{\PP^n\times\PP^n} \PP(\cN_{\ol{M}^1_{\{i\} } | \ol{M}^1}).$
  The fiber of $\ol{M}^2_{\{i\} }\to \ol{M}^1_{\{i\} }$ over
$\{l,q_1,q_2\}$ parametrizes all
choices of lines passing through the second marked point $q_2$ of our
line $l$, whereas the fiber in $Bl_{\Delta}\PP^n\times\PP^n$ over a
point in $\Delta$ parametrizes choices of lines passing through the
third marked point of the curve (since
$Bl_{\Delta}\PP^n\times\PP^n\cong \mbox{ Flag }(1,2,n+1)$).

 The following classes are null in the ring $ B^*(\ol{M}_{0,1} (\PP^n, 2) ;\QQ)$:

$$H^{n+1} \mbox{ , } \frac{ (H + \psi )^{n+1} (H + 2 \psi )^{n+1} }{
 \psi}  \mbox{ , }   T_1T_2 \psi \mbox{ , }   T_i(H +\psi )^{n+1},$$
 $$T_iP(T_i) \mbox{ , }
\left.\sum_{i=1}^2 P(t)\right|^{t=T_i}_{t=0} +\frac{(H+2\psi)^{n+1}}{\psi}.$$

 It is easy to see that the second and fourth classes above are in the ideal generated by the others.

  Let $ S := T_1 +T_2$, $ P := T_1 T_2  $ and $s_{k+1} := T_1^k + T_2^k .$
 From the recurrence relation $ s_{k+1} = s_k S-s_{k-1} P$ one deduces
\bean s_k=\sum_{i=0}^{\lfloor k/2 \rfloor}(-1)^i
\left(\begin{array}{c}k-i\\i\end{array}\right) S^{k-2i}P^i.\eean

 After taking invariants with respect to the action of $S_2$ on
$\{T_1, T_2\}$, we find that
$$A^* (\ol{M}_{0,1} (\PP^n , 2); \QQ ) =\frac{ \QQ[H, \psi ,
S,P]}{\cJ}, $$  where $\cJ$ is generated by $H^{n+1} $, $ P \psi $, $T_1P(T_1)+T_2P(T_2)$ and $P(T_1)+P(T_2)+2 \frac{ (H +
\psi )^{n+1} }{ \psi}- \frac{(H + 2 \psi )^{n+1} }{ \psi}$.
  Using formula (3.10) and the fact that $P\psi=0$, we may write the
last two expressions in terms of $S$ and $P$. Indeed,
  $$P(T_1)+P(T_2)=P(S)+\sum_{k=1}^n\sum_{i=1}^{\lfloor k/2 \rfloor}(-1)^i
\left(\begin{array}{c}k-i\\i\end{array}\right)\left(\begin{array}{c}n+1\\k+1\end{array}\right) S^{k-2i}P^iH^{n-k},$$
  $$T_1P(T_1)+T_2P(T_2)=$$ $$
=SP(S)+\sum_{k=1}^n\sum_{i=1}^{\lfloor k/2 \rfloor}(-1)^i\left(\begin{array}{c}k-i\\i\end{array}\right)\left(\begin{array}{c}n+1\\k\end{array}\right) S^{k-2i}P^iH^{n-k+1}.$$
\end{example}

\section{Appendix}

\begin{remark}
%$H_h^{n+1}=0$
  If $I\not=\emptyset$, the element (4) in the Proposition 3.22 may be
simplified to   \bean \prod_{h\in I\cup J} T_h \cdot P^k_{I,J}(\{T_h\}_{h\in
\cP_{\geq d-k}(I,J)}).\eean   Indeed, in this case $T_hT_{h'}=0$ for all $h,h'\in
\cP_{\geq d-k}(I,J)$ and so, it becomes a nice exercise to see that both the
expression (4) in Proposition 3.22 and expression (4.1) above can be written as
finite sums   $$ \sum_{j}  P^k_{I,J}(\{T_h\}_{h\in \Delta_{h_j}})  $$   where
each $h_j$ is an element of $\cP_{\geq d-k}(I,J)$ and $\Delta_{h_j}=\{h\in
\cP_{\geq d-k}(I,J);  h\supseteq h_j\}$.   By comparing coefficients of monomials in
$T_{h_j}$ one checks that the same terms appear in both sums with the same
coefficients.

   Note that when $I=\{h\}$, relation (4.1)
  tells that  $H_h^{n+1}=0$ on $\ol{M}^k_{Jh}$, where $H_h$ is defined
  as in formula (3.7).

  \end{remark}

\begin{remark}

Let $k< k'<d$ and fix  $h, h'\subset\{1,...,d\}$ such that $|h|=d-k$,
$|h'|=d-k'$. Recall that the preimage of $\ol{M}^{k'}(h')\subset
\ol{M}^{k'}$ in $\ol{M}^{k'}_{h}$ has two distinct components
 $\ol{M}^{k'}_h(h'_{in})$ and
$\ol{M}^{k'}_h(h'_{out})$, according to the
position of the cardinal $(d-k')$-sets in raport to $h$: $h'_{in}\subset h$ and
$h'_{out}\cap h=\emptyset$.  The intersection formula      \bea
T_h[\ol{M}^{k'}(h')]=[\ol{M}^{k'}_h(h'_{in})]+[\ol{M}^{k'}_h(h'_{out})]\eea
is a consequence of  Lemma 3.21 and relations (4) in Proposition 3.22. 
Indeed, by Lemma 3.21 : $$[ \ol{M}^{k'}(h')]=\sum_{h''}
P^{k'}_{\emptyset,\emptyset}(\{T_{h'''}\}_{h''' \supset h''},
t_{h''})\left.\right|^{t_{h''}=T_{h''}}_{t_{h''}=0}+
P^{k'}_{\emptyset,\emptyset}(0),$$ $$
[\ol{M}^{k'}_h(h'_{out})]=T_h(\sum_{h''\not\subseteq h}
P^{k'}_{\emptyset,h}(\{T_{h'''}\}_{h'''\supset h''},
t_{h''})\left.\right|^{t_{h''}=T_{h''}}_{t_{h''}=0}+ P^{k'}_{\emptyset,h}(0)),
$$ $$[\ol{M}^{k'}_h(h'_{in})]=T_h(\sum_{h''\subset h}
P^{k'}_{\emptyset,\emptyset}(\{T_{h'''}\}_{h''' \supset h''},
t_{h''})\left.\right|^{t_{h''}=T_{h''}}_{t_{h''}=0}+
P^{k'}_{\emptyset,\emptyset}(\{T_{h''}\}_{h'' \supseteq h})),$$ where the
sums are taken after $h''$ with $|h''|> d-k'$.  The last term $ T_h
P^{k'}_{\emptyset,\emptyset}(\{T_{h''}\}_{h'' \supseteq h}))$ is the pullback
of the class $[\ol{M}^k_h(h'_{in})]$, as evaluated by the natural
change of variables in formula (3.2): $$H \to H_{h}:=
H+(d-|h|)\psi+\sum_{h''\supset h} |h''\setminus h| T_{h''},$$ $$\psi \to
\psi_{h}(T_h):= \psi+\sum_{h''\supseteq h} T_{h''},$$  and $d \to |h|$ for 
$\ol{M}_{0,1}(\PP^n,|h|)$.

 Whenever $h\cap h''=\emptyset$, the polynomial  $
P^{k'}_{\emptyset,h}(\{T_{h'''}\}_{h''' \supset h''}, t_{h''})$ differs from  $
P^{k'}_{\emptyset,\emptyset}(\{T_{h'''}\}_{h''' \supset h''}, t_{h''})$ by a
multiple of $\psi_{h\cup h''}T_{h\cup h''}$, which becomes zero on
$T_hT_{h''}$ by relation (4) in Proposition 3.22, written for partitions
$(I,J)= $ $(\{h, h''\}, \emptyset)$ at the $(d-|h\cup h''|)$--th intermediate
step. Whenever $h''\supset h$, the difference is summarized by the formula:
\bean \sum_{h''\supset h}  (P^{k'}_{\emptyset,\emptyset}-P^{k'}_{\emptyset,h}
)(\{T_{h'''}\}_{h''' \supset h''},
t_{h''})\left.\right|^{t_{h''}=T_{h''}}_{t_{h''}=0}+ \eean \bea + 
P^{k'}_{\emptyset,\emptyset}(0)-P^{k'}_{\emptyset,h}(0)   =
P^{k'}_{\emptyset,\emptyset} (\{T_{h''}\}_{h'' \supset h},
t_{h})\left.\right|_{t_{h}=0}\eea derived as in the previous remark.
\end{remark}

\begin{lemma}

Let $I\subset \cP$ be a nested set and let $h, h', h''\in \cP\setminus I$
denote elements such that $|h|=d-k$ and $|h'|=|h''|=d-k'$, with $k'< k$.
Assume for simplicity that the $k'$--th and $k$--th steps of blow--up are
performed successively (such that the formulas do not explicitely contain 
exceptional divisors originated in other blow--up steps).  The
following expression  $$ \sum_{h'}
(P_{\ol{M}^{k'}_{Ih'}|\ol{M}^{k'}_{I}}(T_{h'})-P_{\ol{M}^{k'}_{Ih'}|\ol{M}^{k'}_{I}}(0))
 +[\ol{M}^{k'}_{I}(h)]$$ in $B^*(\ol{M}^{k'}_{I})[(T_{h'})_{h'},
(T_{h})_{h}]$  is in the ideal generated by: \begin{enumerate} \item $\sum_{h}
P_{\ol{M}^{k}_{Ih}|\ol{M}^{k}_{I}}(\{T_{h'}\}_{h'\supset h}, t_h)\left.
\right|^{t_h=T_h}_{t_h=0}+ \sum_{h'}
P_{\ol{M}^{k}_{Ih}|\ol{M}^{k}_{I}}(t_{h'})\left.
\right|^{t_{h'}=T_{h'}}_{t_{h'}=0}+[\ol{M}^k_I(h)] ;$ \item $T_h(
\sum_{h'\supset
h}(P_{\ol{M}^{k'}_{Ih'}|\ol{M}^{k'}_{Ih}}(T_{h'})-P_{\ol{M}^{k'}_{Ih'}|\ol{M}^{k'}_{Ih}}(0))+
P_{\ol{M}^{k'}_{Ih'}|\ol{M}^{k'}_{Ih}}(0));$ \item $T_{h'}(\sum_{h; h\cap
h'=\emptyset} P_{\ol{M}^{k}_{Ih}|\ol{M}^{k}_{I}}(\{T_{h''}\}_{h''\supset h},
t_h)\left.  \right|^{t_h=T_h}_{t_h=0}+$

$+\sum_{h''\not=h'}
P_{\ol{M}^{k}_{Ih}|\ol{M}^{k}_{I}}(t_{h''})\left. 
\right|^{t_{h''}=T_{h''}}_{t_{h''}=0}+[\ol{M}^{k}_{Ih'}(h)])$, for $h$
such that $h\cap h'=\emptyset$ ; \item $T_{h'}H_{h'}^{n+1} $ \end{enumerate}

\end{lemma}

 Here $[\ol{M}^k_I(h)]=P_{\ol{M}^{k}_{Ih}|\ol{M}^{k}_{I}}(0)$, where $P_{\ol{M}^{k}_{Ih}|\ol{M}^{k}_{I}}$ is thought of as a function
 of both $t_h$ and $\{t_h''\}_{h''\supset h}$.

 \begin{proof}

%%de regandit in relatiile 1 si 3 cine e cu ~. Rel. 2 nu tocmai, ar trebui polynomul specific lui h'.

  We note that $(1)-(3) = (3')$, where $(3')$ is
\bea T_{h'}(\sum_{h; h\subset
h'}P_{\ol{M}^{k}_{Ih}|\ol{M}^{k}_{I}}(\{T_{h''}\}_{h''\supset h}, t_h)\left.
 \right|^{t_h=T_h}_{t_h=0}+P_{\ol{M}^{k}_{Ih}|\ol{M}^{k}_{I}}(T_{h'}) ), \eea

 and $P_{\ol{M}^{k}_{Ih}|\ol{M}^{k}_{I}}(T_{h'}) = [\ol{M}^{k}_{Ih'}(h)]$, for $h$
such that $h\subset h'$.

 With this, the lemma follows by a direct computation modulo the relations (1)--(4). 
The initial observation is that
  $$P_{\ol{M}^{k'}_{Ih'}|\ol{M}^{k'}_{I}}(t_{h'})\left. \right|^{t_{h'}=T_{h'}}_{t_{h'}=0}= 
P_{\ol{M}^{k'}_{Ih'}|\ol{M}^{k'}_{Ih}}(t_{h'})
  P_{\ol{M}^{k}_{Ih}|\ol{M}^{k}_{I}}(t_{h'})\left.
\right|^{t_{h'}=T_{h'}}_{t_{h'}=0}.$$ ($P_{\ol{M}^{k'}_{Ih'}|\ol{M}^{k'}_{I}}(t_{h'})$ as a
polynomial is not
$P_{\ol{M}^{k'}_{Ih'}|\ol{M}^{k'}_{Ih}}(t_{h'})P_{\ol{M}^{k}_{Ih}|\ol{M}^{k}_{I}}(t_{h'})$ in general, due
to the correction term $H_{h'}^{n+1}$ in the formula for
$P_{\ol{M}^{k'}_{Ih'}|\ol{M}^{k'}_{I}}(t_{h'})$). The right hand side expression above can be
written as \bea P_{\ol{M}^{k'}_{Ih'}|\ol{M}^{k'}_{Ih}}(t_{h'})\left.
\right|^{t_{h'}=T_{h'}}_{t_{h'}=0}P_{\ol{M}^{k}_{Ih}|\ol{M}^{k}_{I}}(T_{h'})+
P_{\ol{M}^{k'}_{Ih'}|\ol{M}^{k'}_{Ih}}(0) P_{\ol{M}^{k}_{Ih}|\ol{M}^{k}_{I}}(t_{h'})\left.
\right|^{t_{h'}=T_{h'}}_{t_{h'}=0}\eea
 Modulo relation $(3')$, this is \bea
& -P_{\ol{M}^{k'}_{Ih'}|\ol{M}^{k'}_{Ih}}(t_{h'})\left.
\right|^{t_{h'}=T_{h'}}_{t_{h'}=0}\sum_{h\subset
h'}P_{\ol{M}^{k}_{Ih}|\ol{M}^{k}_{I}}(\{T_{h''}\}_{h''\supset h}, t_h)\left. 
\right|^{t_h=T_h}_{t_h=0}+&\\ & +P_{\ol{M}^{k'}_{Ih'}|\ol{M}^{k'}_{Ih}}(0)
P_{\ol{M}^{k}_{Ih}|\ol{M}^{k}_{I}}(t_{h'})\left.
\right|^{t_{h'}=T_{h'}}_{t_{h'}=0}.& \eea  
As $T_h$ divides
$P_{\ol{M}^{k}_{Ih}|\ol{M}^{k}_{I}}(\{T_{h''}\}_{h''\supset h}, t_h)\left. 
\right|^{t_h=T_h}_{t_h=0}$, summing the above after all $h'$--s yields,
modulo relation (2):  \bea  & \sum_{h'}
P_{\ol{M}^{k'}_{Ih'}|\ol{M}^{k'}_{I}}(T_{h'})\left.
\right|^{t_{h'}=T_{h'}}_{t_{h'}=0}=&
\\ &=P_{\ol{M}^{k'}_{Ih'}|\ol{M}^{k'}_{Ih}}(0)
\sum_{h}P_{\ol{M}^{k}_{Ih}|\ol{M}^{k}_{I}}(\{T_{h''}\}_{h''\supset h},
t_h)\left.  \right|^{t_h=T_h}_{t_h=0}+& \\ & +
P_{\ol{M}^{k'}_{Ih'}|\ol{M}^{k'}_{Ih}}(0)
P_{\ol{M}^{k}_{Ih}|\ol{M}^{k}_{I}}(t_{h'})\left.
\right|^{t_{h'}=T_{h'}}_{t_{h'}=0},& \eea which, modulo relation (1) equals
$-P_{\ol{M}^{k'}_{Ih'}|\ol{M}^{k'}_{Ih}}(0)[\ol{M}^{k'}_I(h)]=[\ol{M}^{k'}_I(h')]$. \end{proof}

\begin{lemma}
 Let $I$, $J$, be two partial partitions of $M=\{1,...,d\}$, such
that
   $I\cup J$ is also a partial partition.
   With the notations introduced in formula (3.7), the
expressions $$\prod_{h\in I\cup J} T_h P_{I,J}^{d-l_I-1}(\{T_{h'}\}_{h'\in
\cP_{\geq l_I+1}(I, J)}), \mbox{ }  \prod_{h\in I} T_h (\psi +
\sum_{h'\supseteq \bigcup_{h\in I} h} T_{h'})^{|I|-1} $$ ( if $|I|>1)$ and  $
\prod_{h\in I} T_h \psi^{|I|-1}\prod_{i=1}^{d-l_I}(H+i\psi)^{n+1}$  are in the
ideal generated by  \begin{enumerate}
\item $H^{n+1};$
\item $T_hT_{h'}$ if $h\cap h'\not= \emptyset, h$ or $h'$;
\item $T_hT_{h'}(\psi+\sum_{h''\supseteq h\cup h'} T_{h''});$ for
any $h, h'\not=\emptyset$ such that $h\cap h'=\emptyset $;
\item $T_h(\sum_{h'\not\subseteq h}P^{d-1}_{\emptyset, h}(\{T_{h''}\}_{h''\supset h'}, t_{h'})
\left. \right|^{t_{h'}=T_{h'}}_{t_{h'}=0}+P^{d-1}_{\emptyset,
h}(0)).$
 \end{enumerate}

\end{lemma}

\begin{proof}
 The first claim can be reduced to the case when $|I|=1$ and
$J=\emptyset$. Indeed, let
$$\psi_{I,J}:=\psi + \sum_{h'\in \cP_{\geq l_I+1}(I, J)}T_{h'},$$
$$H_{I,J}:= H+(d-\sum_{h\in I\cup J}|h|)\psi+\sum_{h'\in
\cP_{\geq l_I+1}(I, J)}|h'\backslash \bigcup_{h\in I\cup J}h| T_{h'}. $$
Then $P_{I,J}^{d-l_I-1}(\{T_{h'}\}_{h'\in \cP_{\geq l_I+1}(I,J)}) =
\psi_{I,J}^{|I|-1} H_{I,J}^{n+1}$. Assuming that $|I|\geq 2$, we
can choose two distinct $h_1, h_2\in I$. Modulo relation (3),
$$T_{h_1}T_{h_2}\psi_{I,J}=-T_{h_1}T_{h_2}\sum_{h'\supset h_1\cup
h_2, h'\not\in \cP_{\geq l_I+1}(I, J)} T_{h'}.$$ Furthermore, for each $h'$
as above such that $\cup_{h\in I}h\not\subseteq h'$, we can choose
$h_3\in I$ disjoint from $h'$ and, modulo relations (3) and (2),
%$T_{h'}T_{h_3}(\psi+\sum_{h''\supseteq h'\cup h_3} T_{h''})$,
obtain $T_{h'}T_{h_3}\psi_{I,J}=-T_{h'}T_{h_3}\sum_{h''\supset
h'\cup h_3, h''\not\in \cP_{\geq l_I+1}(I,J)} T_{h''}$. Continuing the
same way, we come to the case when $h'\supseteq
\cup_{h\in I}h$. Since $h'\not\in \cP_{\geq l_I+1}(I,J)$, it may be that
either $h'=\cup_{h\in I}h$ or $h'=\cup_{h\in I}h\cup h''$, where
$\{ h\in J;  h\subseteq h''\}\not=\emptyset $, or otherwise $\prod_{h\in J} T_h
T_{h''}=0$. Let $I'=\{h'\}$ and let $ J'=J$ in the first case, $J'=J\setminus
\{h; h\subseteq h'' \}$ in the second case. In both cases,
$T_{h'}H_{I,J}=T_{h'}H_{I',J'}$ and moreover, this also equals
$T_{h'}H_{I',\emptyset}$ after adjusting by
$T_{h'}T_{j'}(\psi+\sum_{h\supseteq h', j'}T_h)$ for every $j'\in J'$.

By the same method, expression $\prod_{h\in I} T_h (\psi +
\sum_{h'\supseteq \cup_{h\in I} h} T_{h'})^{|I|-1} $ is reduced to
relation (3).

 Consider now $I=\{h\}$, $J=\emptyset$. Let
  $\psi_h(t_h):=\sum_{h'\supset h} T_{h'}+t_h$. We proceed by decreasing induction on $|h|$.
  Noting that in all our expressions, $\psi$ plays the role of $T_h$ with $h=\{1,...,d\}$, we can start the induction from the relation $H^{n+1}=0$.
   We can now work out the general induction step.
   Multiplication of
  expression (4) by $\psi_h(0)$ gives, modulo expressions (2) and
  (3):
  \bea &T_h (\sum_{h'\cap h=\emptyset} P^{d-1}_{\emptyset,h}(\{T_{h''}\}_{h''\supset h'}, t_{h'})
\left. \right|^{t_{h'}=T_{h'}}_{t_{h'}=0}(\sum_{h''\supset h,
h''\cap h'=\emptyset}T_{h''})+  &  \\ & + \sum_{h'\supset h}
P^{d-1}_{\emptyset,h}(\{T_{h''}\}_{h''\supset h'}, t_{h'}) \left.
\right|^{t_{h'}=T_{h'}}_{t_{h'}=0}(\sum_{h'\supset h''\supset
h}T_{h''}+\psi_{h'}(T_{h'}))+&\\ & +P^{d-1}_{\emptyset, h}(0)\psi_h(0)
),&\eea which, after regrouping of terms becomes \bean & T_h (
\sum_{h''\supset h}T_{h''}\sum_{h'\not\subseteq h''}
P^{d-1}_{\emptyset,h}(\{T_{h'''}\}_{h'''\supset h'}, t_{h'}) \left.
\right|^{t_{h'}=T_{h'}}_{t_{h'}=0}+& \\ & + \sum_{h'\supset h}
\psi_{h'}(T_{h'})P^{d-1}_{\emptyset,h}(\{T_{h''}\}_{h''\supset h'},
t_{h'}) \left. \right|^{t_{h'}=T_{h'}}_{t_{h'}=0}
+P^{d-1}_{\emptyset, h}(0)\psi_h(0) ).& \eean
 As remarked earlier, if $h'\cap h''=\emptyset$, then $$P^{d-1}_{\emptyset,h}(\{T_{h'''}\}_{h'''\supset h'}, t_{h'})
\left.
\right|^{t_{h'}=T_{h'}}_{t_{h'}=0}=P^{d-1}_{\emptyset,h''}(\{T_{h'''}\}_{h'''\supset
h'}, t_{h'}) \left. \right|^{t_{h'}=T_{h'}}_{t_{h'}=0}$$ and by formula
(4.2)\bea & \sum_{h'\supset
h''}(P^{d-1}_{\emptyset,h}-P^{d-1}_{\emptyset,h''})(\{T_{h'''}\}_{h'''\supset
h'}, t_{h'}) \left.
\right|^{t_{h'}=T_{h'}}_{t_{h'}=0}+P^{d-1}_{\emptyset,h}(0)-P^{d-1}_{\emptyset,h''}(0)= &\\
& =P^{d-1}_{\emptyset,h}(\{T_{h'''}\}_{h'''\supset h'}, t_{h''})
\left. \right|_{t_{h''}=0} & \eea modulo (2) and (3) and thus,
modulo expression (4) written for $h''$, the summand (4.3) above
is \bea -T_h\sum_{h''\supset
h}T_{h''}(P^{d-1}_{\emptyset,h}(0)+P^{d-1}_{\emptyset,h}(\{T_{h'''}\}_{h'''\supset
h'}, t_{h''})\left. \right|_{t_{h''}=0}). \eea On the other hand, summand (4.4)
may be written as \bea & T_h ( \sum_{h'\supset h}
\psi_{h'}P^{d-1}_{\emptyset,h}(\{T_{h''}\}_{h''\supset h'}, t_{h'})
\left. \right|^{t_{h'}=T_{h'}}_{t_{h'}=0}+ &\\ &+  \sum_{h'\supset h}
\psi_{h'}P^{d-1}_{\emptyset,h}(\{T_{h''}\}_{h''\supset h'}, t_{h'})
\left. \right|_{t_{h'}=0}(\psi_{h'}(T_{h'})-\psi_{h'}(0))+&\\ &
+P^{d-1}_{\emptyset, h}(0)\psi_h(0) ).& \eea As
$\psi_{h'}(T_{h'})-\psi_{h'}(0)=T_{h'}$, putting the two summands
together we obtain \bea \sum_{h'\supset h}
\psi_{h'}P^{d-1}_{\emptyset,h}(\{T_{h''}\}_{h''\supset h'}, t_{h'})
\left. \right|_{t_{h'}=0}+P^{d-1}_{\emptyset, h}(0)\psi. \eea In
view of definition of $P^{d-1}_{\emptyset, h}$ (formula (3.7)) and of Remark
4.1, the above reduces to $T_hH_{h,\emptyset}^{n+1}$. The
induction hypothesis was necessary when applying Remark 4.1. This ends the
proof for the second expression in the lemma.

 Dependence of the last expression on terms (1)--(4) can be deduced by decreasing induction on $I$:
\bea  &  r(I):=\prod_{h\in I} T_h \psi^{|I|-1}\prod_{i=1}^{d-l_I}(H+i\psi)^{n+1}= \prod_{h\in I} 
T_h P_{\ol{M}^1_I|\ol{M}^1}(\psi) = &\\ &= \prod_{h\in I\cup \{h'\}} T_h
[\ol{M}^1_{Ih'}] P_{\ol{M}^1_{Ih'}|\ol{M}^1}(\psi),& \eea  for any $h'$ such
that $h'\cap (\cup_{h\in I}h)=\emptyset$, where $[\ol{M}^1_{Ih'}]$ is the
class of $\ol{M}^1_{Ih'}$ in $\ol{M}^1_I$.  Thus by relation (4) of
Proposition 3.22, $r(I)$ is in the ideal generated by $r(Ih'),$ $h'$ as
above. The last induction step coincides with the second expression in the
lemma when $I$ is a complete partition of $M$.

%  Let $h\subset \{1,...,d\}$ and let $H_h:=
% H+(d-|h|)\psi+\sum_{h'\supset h}|h'\backslash h| T_{h'}. $ Let
% $\psi_h(t_h):=\sum_{h'\supset h} T_{h'}+t_h$ and let $\psi_h$ be
% the short notation for $\psi_h(T_h)$. As in relation (4.15), we
% denote by $P^{d-1}_{\emptyset, \{h\}}(\{T_{h''}\}_{h''\supset h'},
% t_{h'})$, or, for short, $P^{d-1}_h(t_{h'})$, the polynomial
% $$P^{d-1}_h(t_{h'})=\psi_{h'}(t_{h'})^{-1}((H_{h'}+|h'\backslash
% h|\psi_{h'}(t_{h'}))^{n+1}-H_{h'}^{n+1})$$ and
% $P^{d-1}_{Ih}(0)=\psi^{-1}((H_+(d-|h|)\psi )^{n+1} -H^{n+1})$.
 %  Then the expression $T_hH_h^{n+1}$

\end{proof}

   %References

\providecommand{\bysame}{\leavevmode\hbox to3em{\hrulefill}\thinspace}

\end{document}